\documentclass{cmslatex}
\usepackage[paperwidth=7in, paperheight=10in, margin=.875in]{geometry}
 \usepackage[backref,colorlinks,linkcolor=red,anchorcolor=green,citecolor=blue]{hyperref}
\usepackage{amsfonts,amssymb}
\usepackage{amsmath}
\usepackage{graphicx}
\usepackage{cite}
\usepackage{enumerate}
\usepackage{caption}
\usepackage{subfigure}

\renewcommand{\theequation}{\arabic{section}.\arabic{equation}}

   \allowdisplaybreaks
\begin{document}
 \title{A hybrid numerical method for a microscopic and macroscopic traffic flow model}%\thanks{Received 5 November 2019, accepted 21 December 2020.}}

          %For each author, make a block with the following macros:

\author{Yuanhong Wu\thanks{Department of mathematics and statistics, Ningbo University, 315211, P.R. China (2311400047@nbu.edu.cn).}
\and Shuzhi Liu\thanks{ School of Statistics and Data Science, Ningbo University of Technology, 315211, P.R. China (szliu@nbut.edu.cn).}
	\and Qinglong Zhang\thanks{Department of mathematics and statistics, Ningbo University, 315211, P.R. China, (zhangqinglong@nbu.edu.cn).}
	}
          \
 %         \thanks{address, zhangqinglong@nbu.edu.cn.}
          
         \pagestyle{myheadings} \markboth{traffic flow model}{Y. Wu, S. Liu, and Q. Zhang} \maketitle

          \begin{abstract}
In this paper, we introduce a traffic flow model based on a microscopic follow-the-leader model, while enforcing maximal constraints on the density and velocity of the flow. The related macroscopic model can be represented in conservative formulation. By introducing an advected variable $\tilde{u}p$ with the flow, where $p$ is the velocity offset, and $\tilde{u}$ is the `` relative'' velocity, we reformulate the classical Aw–Rascle-Zhang (ARZ) model \cite{AR_model_2000,Zhang} and the modified Aw-Rascle model \cite{Berthelin_2008} to describe a realistic fundamental diagrams. The elementary waves are derived, and the Riemann problem is solved to validate the model's theoretical consistency. We further extend to a two-dimensional model. Numerical simulations are given for both one- and two-dimensional case by using the hybrid Godunov-Glimm scheme \cite{Betancourt, Christophe_2007_Godunov_Gimm} to verify the model's performance. 
          \end{abstract}


  \section{Introduction}
Modeling traffic flow dynamics accurately is a challenge in the transportation engineering and applied mathematics. As urban mobility systems become increasingly complex, there is a growing demand for reliable traffic models. Historically, a wide range of traffic flow models have been developed, which can be broadly categorized into {\em microscopic}, {\em cellular} and {\em macroscopic} frameworks. {\em Microscopic} models simulate the individual behavior of vehicles by ODE, typical examples are known as ``follow-the-leader" models \cite{Gazis_1961, Bando_1995_OVM, Jiang_Rui_2001, Helbing_1998, Treiber_2013, Gipps_1981, Matin_2023_OVM_FTL, Guang_han_2013_AOVM}. Recent developments on microscopic models which integrate data-driven approaches, including neural networks \cite{Kuefler_2017, Linsen_2013, Mofan_2017}, recurrent architectures \cite{Ma_2015}, and reinforcement learning \cite{Emamifar_2023, Meixin_2020, Meixin_2018}. {\em Cellular} models divide the road into cells and evaluate how vehicles advance through cells \cite{Daganzo2}. {\em Macroscopic} models describe traffic flow from a continuum perspective, treating the aggregation of vehicles as a compressible fluid governed by conservation laws. These models are described by PDEs involving traffic density, velocity, and flux. Examples are listed as the Lighthill-Whitham model \cite{Lighthill_Whitham_1955}, the Payne model \cite{Payne}, and the ARZ model \cite{AR_model_2000, Zhang}. Other familiar approaches are kinetic models (such as the cellular automata models \cite{Klar}). The macroscopic models are based on the balanced laws on the observed equations. In the second order models, there are basically two conservation equations, one is the conservative of total cars, the other is the acceleration equation, although it can be written in a conservative form, but can not be interpreted as the conservation of ``momentum". 
One of the representative model is the ARZ traffic flow model. In \cite{Aw_2002}, the authors proven that the ARZ model is equivalent to a microscopic car-following model via scaling. Through this model effectively resolves several issues identified by Daganzo \cite{Daganzo1}, it does not guarantee that the upper limits for density and velocity are respected during evolution. In \cite{Berthelin_2008}, the authors proposed a modified AR (MAR) model with the maximal density constraint, they showed that as the density approaches the maximum value $\rho^{*}$, the pseudo-velocity $p(\rho)$ tends to infinity. The MAR model is capable of capturing the emergence and evolution of traffic jams.
 
 In this paper, our motivation is to further put the velocity constraint into the MAR model, i.e., the velocity has the upper bound $u^{*}$, and the velocity $u\to u^{*}$ as the density $\rho\to 0$. This can be obtained by introducing a new variable, which can be called the ``pseudo-velocity" variable $\tilde{u} = (\frac{1}{u}-\frac{1}{u^{*}})^{-1}$, for simplicity. However, another problem will be caused by using $\tilde{u}$ instead of $u$ in the MAR model, that is $\tilde{u}+p$ will no longer keep invariant along the traffic trajectories. For example, with the definition of $\tilde{u}$, as $\rho \to 0$, we have $p \to 0$ and $\tilde{u} \to \infty$, hence $\tilde{u} + p \to \infty$. This prevents the smooth transition from the ``vacuum" state to the free flow where $\tilde{u}+p$ is some finite number.
 
 The two limiting cases that $\rho\to 0$, $\tilde{u}\to \infty$, and $\rho\to \rho^*$, $p\to \infty$ indicate that the ARZ model can not be used directly with the new variable $\tilde{u}$ introduced. To cure this deficiency, we begin with a microscopic car-following framework, particularly, we modified the  Gazis-Herman-Rothery (GHR) model \cite{Gazis_1961} by refining the acceleration equation. Then, we go back to the macroscopic model and conclude that $\tilde{u}p$ will be conserved and advected to the traffic trajectories.
  
In complex traffic environments such as urban networks and multilane highways, classical one-dimensional traffic flow models are not enough since traffic behaviors like lane changing, merging, and diverging play an important role. To better understanding such dynamics, two-dimensional macroscopic models have been proposed \cite{Sukhinova_2009, Nagatani_1999}. More recently, Herty, Moutari and Visconti \cite{Herty_2018_2D} proposed a two-dimensional macroscopic model based on the ARZ framework, which was developed by applying temporal scaling techniques to a microscopic car-following system. Indicated by their works, we extend the proposed model to the 2-D scenrio. Numerical simulations show that the 2-D model can simulate 2-D traffic waves and the traffic behaviors such as  lane changing, merging  through wave interactions.
  
The development of robust and efficient numerical methods becomes essential for simulating complex traffic phenomena with high fidelity. Traditional numerical schemes are available such as finite-difference and finite-volume discretizations \cite{LeVeque_2002, Morton_2007}. Among them, the Godunov scheme is attracting since it takes use of the  Riemann solution to construct the numerical flux. However, Chalons and Goatin  \cite{Christophe_2007_Godunov_Gimm} have pointed out that nonphysical oscillations occur near contact discontinuities when the Godunov scheme is applied to the ARZ model, they proposed a hybrid method which combine the random choice method with the Godunov scheme, their method is proven to be quite efficient in capturing contact discontinuities. We apply their hybrid method to the present model, which contain two steps: first, use the Glimm random choice method to resolve the contact discontinuities, and then use the Godunov scheme to resolve nonlinear waves.  Moreover, we extend the hybrid method to the two-dimensional system by using the splitting techniques, and compare the results with that obtained by the HLL scheme. It shows that the hybrid method is efficient in capturing the position of contact discontinuities and the swirling of their interactions.

The paper begins with the model derivation in Section~\ref{sec1}, integrating microscopic and macroscopic viewpoints. By modifying the classical ARZ model and the GHR model, we refine the model by imposing both the maximum density and velocity constraints. We also give comparative experiments with other models, which demonstrate that the revised auxiliary variable yields notable improvements.
In Section~\ref{sec2}, we conduct a detailed characteristic analysis of the proposed model. The elementary waves, namely rarefaction wave, shock wave, and the contact discontinuity are introduced and discussed.
Section~\ref{sec3} examines three numerical schemes: the HLL method, the Godunov method, and a combined Godunov–Glimm approach. Through comparative analysis and simulations, we show that the hybrid method more effectively captures the wave dynamics and aligns well with theoretical expectations.
In Section~\ref{sec4}, we extend the model to two dimensions to account for both longitudinal and lateral traffic waves. The corresponding numerical simulations further illustrate the ability of the model to capture realistic traffic dynamics in multi-lane or urban scenarios.

\section{Model Formulation}\label{sec1}
\subsection{Microscopic Formulation}

Microscopic follow-the-leader model describes vehicles as point particles indexed by $i \in \mathbb{Z} $, each characterized by a time-dependent position $x_i(t)$ and velocity $ u_i(t)$. This model encapsulates a core behavioral assumption: drivers adapt acceleration in response to the relative speed and spacing with respect to the vehicle immediately ahead. The interaction between vehicles is described by the following ordinary differential equation \cite{Chandler_1958}:
\begin{align} \label{FTL}
	\dot{u}_i = C \frac{u_{i+1} - u_i}{(x_{i+1} - x_i)^{\gamma+1}},
\end{align}
where $C>0$ denotes a sensitivity coefficient, and $\gamma > 0 $ determines the influence of headway on driver response. As demonstrated in \cite{Aw_2002}, starting from the microscopic system in Eq.~\eqref{FTL}, the macroscopic ARZ model can be obtained as its hydrodynamic limit. Moreover, the Gazis–Herman–Rothery (GHR) model can be reformulated in the following canonical structure \cite{Gazis_1961}:

\begin{align} \label{GHR}
	\dot{u}_i = \alpha u_i^{\beta} \frac{u_{i+1} - u_i}{(x_{i+1} - x_i)^{\gamma}},
\end{align}
where $\alpha, \beta, \gamma$ are calibration parameters. However, this model neglects the variable sensitivity of drivers to different levels of speed differences. When the ratio $(u_{i+1} - u_i)/{(x_{i+1} - x_i)^{\gamma+1}}$ is fixed at a constant level, the model ``relaxes" to an acceleration equation, with the right side of which proportional to $\alpha u_i^\beta$. Thus the acceleration equation increases monotonically with the follower's velocity $u_i$, such uniform response is inconsistent with empirical driving behavior \cite{Newell_1961}.
 
To address this limitation, we incorporate a quadratic acceleration regulation term. The refined acceleration law is given by:
\begin{align}\label{mar}
	\dot{u}_i = \frac{\gamma u_i(u^* - u_i)}{u^*} \frac{u_{i+1} - u_i}{x_{i+1} - x_i- d}.
\end{align}
where $\gamma$ is the calibration coefficient, $ u^*$ represents the maximum speed and $d$ denotes the minimal admissible inter-vehicle distance, respectively. To compare the model with others, we abbreviate it as RARZ (refined ARZ) hereafter.
The term $x_{i+1} - x_i- d$ in the denominator captures the effective spacing while inducing a repulsive effect as vehicles approach the minimal admissible distance.
In the traffic congestion as vehicle speed approaches zero ($u_i \to 0$), $\gamma u_i(u^* - u_i)/u^* \approx \gamma u_i$ as the quadratic term can be neglected. Then we recover the GHR model with the unit parameter $\beta$. 
 In the free flow as the density tends to zero, the acceleration term will asymptotically approaches zero since the vehicle reaches the maximum allowed speed ($u_i \to u^*$), which ensure the convergence to the maximum speed without overfitting. The modification in the right side of \eqref{mar} has essential significance 
 when the velocity difference $u_{i+1} - u_i$ is proportional to the distance $x_{i+1} - x_i - d$, 
 %when the ratio $(u_{i+1} - u_i)/(x_{i+1} - x_i - d)$ remains constant,
 it reveals that the acceleration regulation term $\gamma u_i(u^* - u_i)/u^*$ reaches its maximum at $u_i = u^*/2$, endowing the model with velocity adjustment sensitivity during traffic phase transitions (e.g., between free-flow and congested states). Another advantage of the model is in the fundamental diagram, it can reveal the bivariate equilibrium relations of the fundamental diagram under varying maximum speed settings. The model thus inherently adapts to the complex dynamical constraints across traffic flows while preserving stability and flexibility, as we will demonstrate in the following part. 
 
\subsection{Macroscopic Formulation}
This section is devoted to deriving the macroscopic model of \eqref{mar} through a scaling limit process.
From a microscopic viewpoint, for any vehicle $i$ with trajectory $x_i(t)$ in a stream of identical particles (length $d$), the local microscopic density is:
$$\rho_i = \frac{\Delta X}{x_{i+1} - x_i},$$
where $\Delta X$ denotes the effective vehicle length that accounts for safety margins. Under bumper-to-bumper condition$(x_{i+1} - x_i = d )$, this density attains its maximum value $\rho^*=\Delta X/d$. Then the normalized relative distance is expressed as:
$$\tau_i = \frac{x_{i+1} - x_i -d}{\Delta X} = \frac{1}{\rho_i} - \frac{1}{\rho^*}.$$

Assuming continuous differentiability of trajectories $x_i(t)$, the temporal evolution of $\tau_i$ can be expressed through a finite difference relation:
$$ \frac{\tau_i(t+\Delta t) - \tau_i(t)}{\Delta t} = \frac{1}{\Delta X} \left( \frac{ x_{i+1}(t+\Delta t) - x_{i+1}(t)}{\Delta t} - \frac{x_i(t+\Delta t) - x_i(t)}{\Delta t} \right).$$

In the limit of vanishing time step $\Delta t \to 0$, this discrete relation recovers the continuous form:
$$\dot{\tau}_i = \frac{\dot{x}_{i+1} - \dot{x}_i}{\Delta X} = \frac{u_{i+1} - u_i}{\Delta X},$$
where $\dot{\tau}_i$ denotes the time derivative of $\tau_i$, and $u_i = \dot{x}_i$ represents the velocity of the $i$-th vehicle.

Let us define
$$w_i=\widetilde{u}_i p(\tau_i), \quad {\rm where} \  \widetilde{u}_i = \left(\frac{1}{u_i} - \frac{1}{u^*}\right)^{-1}, p(\tau_i)= \tau_i^{-\gamma} = \left(\frac{1}{\rho_i} - \frac{1}{\rho^*}\right)^{-\gamma},$$
and using the refined acceleration law given in \eqref{mar}, the time derivative of $w_i$ becomes
\begin{align*}
	\dot{w}_i =& (\widetilde{u}_i p(\tau_i))' \\
	=& u_i^{-2} \widetilde{u}_i^2 p(\tau_i) \dot{u}_i - \gamma \widetilde{u}_i \tau_i^{-1} p(\tau_i) \dot{\tau}_i \\
	=& \widetilde{u}_i p(\tau_i)\left( u_i^{-2}\widetilde{u}_i\dot{u}_i - \gamma \tau_i^{-1} \dot{\tau}_i \right) \\
	=& \widetilde{u}_i p(\tau_i)\left( \frac{u^*}{u_i(u^*-u_i)} \frac{\gamma u_i(u^* - u_i)(u_{i+1} - u_i)}{u^*(x_{i+1} - x_i -d)} - \gamma \frac{\Delta X}{ x_{i+1} - x_i - d } \frac{u_{i+1} - u_i}{\Delta X} \right)\\
	=&0.
\end{align*}

Therefore, $\dot{w}_i =0$, which implies that $w_i=\widetilde{u}_i p(\tau_i)$ is conserved in time. By taking the continuum limit and applying the standard transformation from microscopic to macroscopic quantities, we arrive at the following macroscopic model
\begin{align} \label{model1}
	\left\{
	\begin{aligned}
		& \rho_t +(\rho u)_x=0,\\
		& (\rho \widetilde{u} p)_t+(\rho u \widetilde{u} p)_x=0.
	\end{aligned}
	\right. 
\end{align}
where $ \widetilde{u} = \left(\frac{1}{u} - \frac{1}{u^*}\right)^{-1}, p=\left(\frac1\rho-\frac{1}{\rho^*}\right)^{-\gamma}$.

Let's compare the present model (abbreviated as RARZ) with the ARZ model and MAR model from the macroscopic point of view before any further analysis. For completeness, we post the other two model. The ARZ model \cite{AR_model_2000, Zhang} is formulated as:
\begin{align} \label{AR}
	\left\{
	\begin{aligned}
		& \rho_t +(\rho u)_x=0,\\
		& (\rho w)_t+(\rho u w)_x=0,
	\end{aligned}
	\right. 
\end{align}  
where the variable $ \omega= u + p(\rho)$ is advected with the flow and travels along the vehicle's trajectories.  
However, this standard formulation does not adequately incorporate physical constraints of maximum speed $u^*$ and maximum density $\rho^*$. To fix the deficiency, Berthelin et al. \cite{Berthelin_2008} proposed a modified pressure function of the form:
\begin{align}\label{p}
	p(\rho) = \left( \frac{1}{\rho} - \frac{1}{\rho^*}\right)^{-\gamma},
\end{align}  
which introduces a singularity of $p(\rho)$ as $\rho \to \rho^*$. The modified term effectively keep the upper bound of the density, which will not exceed the congestion density.

Based on the MAR model, the RARZ model present here further introduces the maximum allowed speed $u^* $, and introduce a ``relative"  velocity defined by: 
\begin{align}\label{u}
	\widetilde{u} = \left( \frac{1}{u} - \frac{1}{u^*} \right)^{-1}.
\end{align}
Clearly, $\widetilde{u}$ is monotonically increasing, which satisfies $u\to 0$, $\widetilde{u} \to 0$, and $u\to u^*$, $\widetilde{u} \to \infty$. The introduction of $\widetilde{u}$ is not trivial  since we can not just replace $u$ with $\widetilde{u}$ in the MAR model. The advected variable $ \widetilde{u} + p(\rho)$ is not expected to keep unchanged apart from the contact discontinuity. Otherwise, take $\rho\to 0$, then $u$ tends to $u^{*}$, $\widetilde{u} \to \infty$, and thus $\widetilde{u} + p(\rho) \to \infty$, which means that even in the free flow, across the rarefaction wave, $\widetilde{u} + p(\rho)$ will not remain constant. 
 
In the present model \eqref{model1}, the variable  $\widetilde{w} = \widetilde{u} p(\rho)$ is convected to the traffic flow where both density and velocity constraints are integrated into $\widetilde{w}$. Such formulation answers the above concern that in the free flow as $\rho\to 0$, $p(\rho)\to 0$, and $\widetilde{u} \to \infty$, then the limit of $\widetilde{u} p(\rho)$ has a indeterminate form. Simultaneously, in the congested regimes where the density approaches maximal value $\rho^*$, $p(\rho) = \left( \frac{1}{\rho} - \frac{1}{\rho^*} \right)^{-\gamma} \to \infty$, and $\widetilde{u} \to 0$. The limit of $\widetilde{u} p(\rho)$ can still remain to a finite number if the parameters are properly chosen. This enables a smooth transition between free-flow and congested regimes.

\begin{figure}[h!]
  	\vspace{-0.3cm}
 	\centering
 	\includegraphics[width=1.0\linewidth]{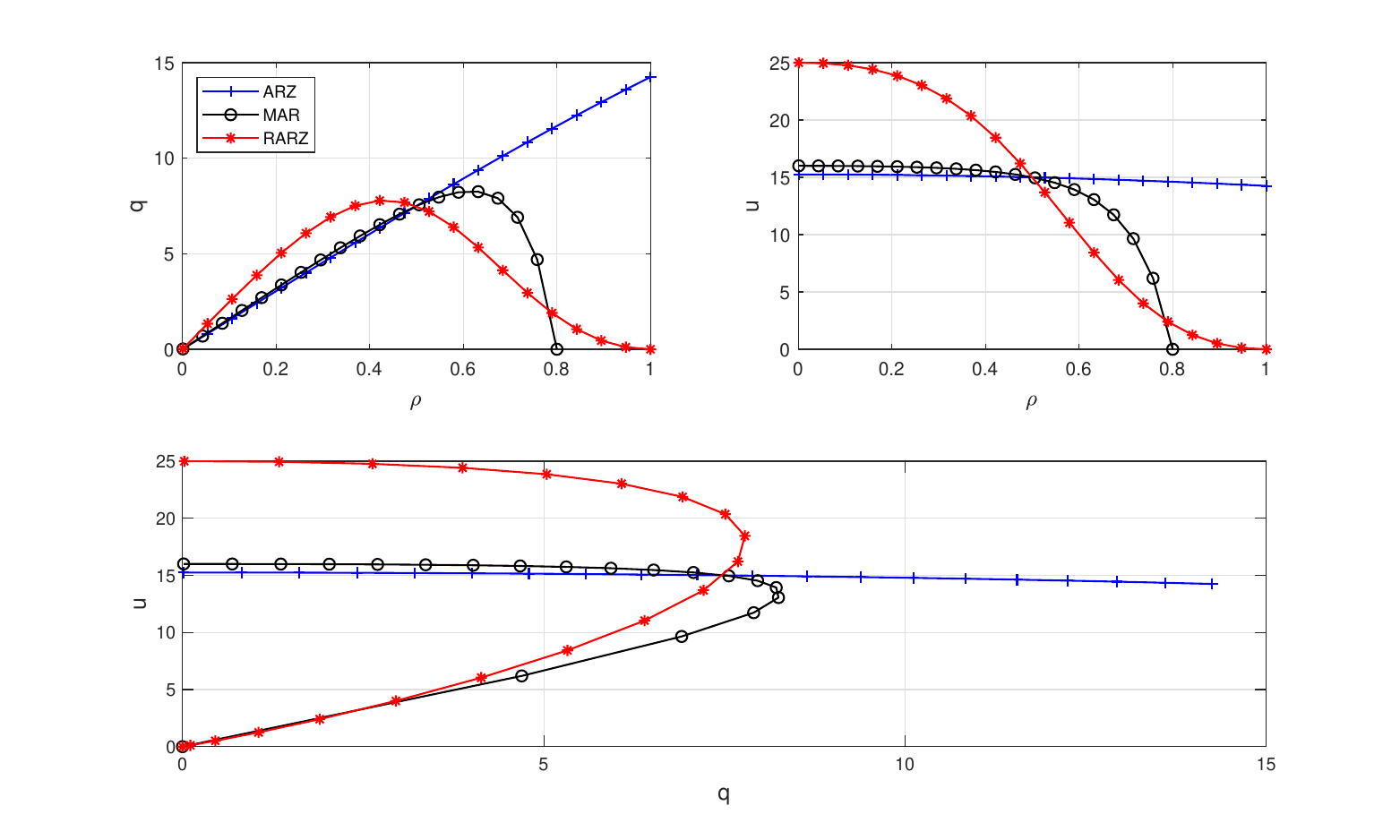}
 	\vspace{-0.8cm}
 	\caption{Comparison of fundamental diagrams among the three models.}
 	\label{fdmc}
 	\vspace{-0.9cm}
 \end{figure}

A comparative analysis of the fundamental diagrams for the ARZ, MAR, and RARZ models is presented below. The fundamental diagrams connect with traffic density, velocity and flow through an empirical or theoretical perspective \cite{Lighthill_Whitham_1955, Richards_1956, Treiber_2013}. Fig.~\ref{fdmc} present the fundamental diagrams from theoretical analysis. For the ARZ model, The velocity offset is defined by $p(\rho)=\rho^\gamma$, leading to the conserved quantity $w=u+p(\rho)$. However, under high-density conditions, this formulation fails to enforce a realistic velocity decay, which result in relative high speeds despite in the congestion traffic flow,  thus violate the physical constraints. The MAR model replaces the velocity offset to \eqref{p}, which properly diverges as $\rho \to \rho^*$, mimicking the bumper-to-bumper condition. However, the conserved quantity $w=u+p(\rho)$ in this setting may lead to $u \to - \infty$, in order to maintain the conservation of $u+p$ across the shock wave or the rarefaction wave. Thus, negative speeds will be produced, which are nonphysical. The RARZ model resolves both the density and the velocity constraints \eqref{mar}, ensuring non-negative speeds and a physically realistic response under congested conditions, as shown in Fig.~\ref{fdmc}.

   \begin{figure}[h!]
 	\centering
 	\includegraphics[width=1.0\linewidth]{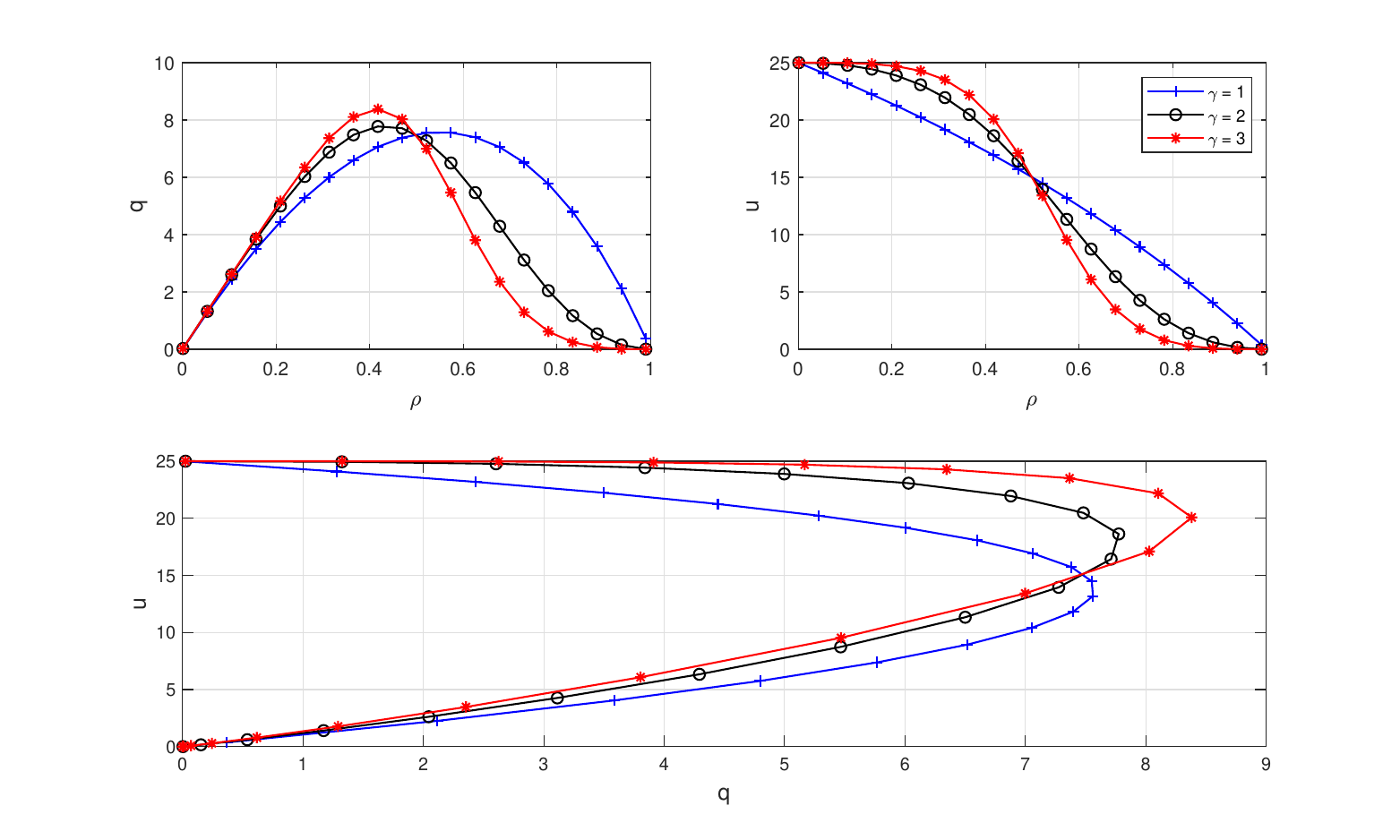}
 	\vspace{-0.8cm}
 	\caption{Comparison of fundamental diagrams with different $\gamma$.}
 	\label{fdgamma3}
 	\vspace{-0.5cm}
 \end{figure}
 
Figure~\ref{fdgamma3} present an analysis of the fundamental diagram of RARZ model under different values of the parameter $\gamma$. As $\gamma$ increases, the velocity decreases more sharply with respect to the increasing of density, which represent a more strongly reaction to traffic congestion. This represents the model's flexibility in capturing a wide range of driver behaviors and traffic conditions. Moreover, the first subfigure in \eqref{fdgamma3} visually resembles an inverted $\lambda$ shape\cite{Siebel}, especially when $\gamma = 3$.

\section{Characteristic analysis}\label{sec2}
\subsection{Preliminaries}
Introducing the vector $U=(\rho, u)^{T}$, system \eqref{model1} can be rewritten compactly as:
$$ A(U)\partial_t U + B(U)\partial_x U=0,$$ 
where
\begin{align*}
	A(U)=
	\begin{pmatrix}
		1 & 0 \\
		\frac{(\gamma+1)\rho^* - \rho}{\rho^* - \rho} \widetilde{u} p & \rho p \left( \frac{u^*}{u^* - u}\right)^2
	\end{pmatrix}, 
	 \ B(U)=
	\begin{pmatrix}
		u & \rho \\
		\frac{(\gamma+1)\rho^* - \rho}{\rho^* - \rho} u \widetilde{u} p &  \frac{2u^* - u}{u^* - u}\rho \widetilde{u} p
	\end{pmatrix}.
\end{align*}
The system has two eigenvalues:
$$ \lambda_1=u , \quad \lambda_2=u - \frac{\gamma u \rho^* \left(u^* - u\right)}{u^*\left(\rho^* -\rho \right)}.$$
This reveals the inherent anisotropy of traffic flow, in which information cannot propagate faster than the driver's speed. 
The right eigenvectors related to this system are given by:
$$ \vec{r}_1=(1, \ 0)^T, \quad \vec{r}_2=\left(\frac{\rho u^* \left(\rho^* - \rho\right)}{\gamma u \rho^* \left(u^* - u\right)}, \ -1 \right)^T.$$

Since $\nabla \lambda_1 \cdot \vec{r}_1 = 0$ and $\nabla \lambda_2 \cdot \vec{r}_2 \neq 0$, it follows that the first characteristic field is linearly degenerate, while the second characteristic field exhibits genuine nonlinearity. System \eqref{model1} is strictly hyperbolic unless $u=u^*$,  which means the two  eigenvalues coincide when the velocity reaches the maximum constraint speed.

\subsection{The rarefaction waves}
For the self-similar rarefaction wave, which depends on $\xi = x/t$, the Riemann invariants $w_i$ associated with the eigenvalues $\lambda_i$ are computed by solving
$$\nabla w_i \cdot \vec{r}_i = 0, \quad i = 1, 2.$$

For $ \lambda_1$, the corresponding Riemann invariant is given by
$$ w_1 = u. $$
For $ \lambda_2 $, the corresponding Riemann invariant is given by
$$ w_2 = ln \frac{u}{u^* - u} + \gamma \frac{\rho}{\rho^* - \rho}. $$
As $\lambda_1$ is linearly degenerate, thus the rarefaction wave is linked to the $\lambda_2$ characteristic. Note that $ w_2 =const. $ is equivalent to $\widetilde{u} p = const$. Therefore, for a specified left-hand state $U_0 = (u_0, \rho_0)$, the corresponding right-hand state $U = (u, \rho)$ that can be reached via the rarefaction curve is

\begin{align}\label{R}
	R(U, U_0):
	\left\{
	\begin{aligned}
		& \xi = \lambda_1 = u - \frac{\gamma u \rho^* (u^* - u)}{u^*(\rho^* - \rho)}, \\
		& \widetilde{u} p = \widetilde{u}_0 p_0, \\
		& \rho < \rho_0, \quad u > u_0.
	\end{aligned}
	\right.
\end{align}

\subsection{The discontinuous waves}
Considering a discontinuity traveling at speed ${\rm d}x/{\rm d}t = \sigma$, the Rankine–Hugoniot condition for system~(\ref{model1}) takes the form
\begin{align}\label{RH}
	\left\{
	\begin{aligned}
		& \sigma[\rho]=[\rho u] ,\\
		& \sigma[\rho \widetilde{u} p]=[\rho u \widetilde{u} p].
	\end{aligned}
	\right.
\end{align}
The notation $[f] := f_1 - f_0$ is used to indicate the jump of $f$ across the discontinuity, while $\sigma$ corresponds to the wave speed linking the states $U_0 = (\rho_0, u_0)$ and $U_1 = (\rho_1, u_1)$.

By solving (\ref{RH}), one has two types of discontinuities solutions, the contact discontinuity and the shock wave. 

Starting from a left state $U_0 = (\rho_0, u_0)$, the right state $U = (\rho, u)$ associated with a contact discontinuity satisfies
$$J(U, U_0): \quad u = u_0 = \xi, \quad \rho \ne \rho_0.$$

The shock wave connecting $U_0$ and $U$ satisfies the following conditions:
\begin{align}\label{S}
	S(U, U_0): \quad
	\left\{
	\begin{aligned}
		& \sigma = \frac{\rho u - \rho_0 u_0}{\rho - \rho_0}, \\
		& \widetilde{u} p = \widetilde{u}_0 p_0, \\
		& \rho > \rho_0, \quad u < u_0,
	\end{aligned}
	\right.
\end{align}
here, $\sigma$ denotes the shock speed, while the second relation ensures the invariance of the variable $\widetilde{u} p$ across the shock. Moreover, for the shock wave, we have 
$$\left(\frac{1}{u} - \frac{1}{u^*}\right)^{-1}= \widetilde{u}_0 p_0 \left(\frac{1}{\rho} - \frac{1}{\rho^*}\right)^{\gamma}.$$
By differentiating it on both sides, we get
$$ \frac{1}{u^2} \left(\frac{1}{u} - \frac{1}{u^*}\right)^{-2} \mathrm{d} u = - \frac{\gamma \widetilde{u}_0 p_0}{\rho^2} \left(\frac{1}{\rho} - \frac{1}{\rho^*}\right)^{\gamma-1} \mathrm{d} \rho,$$
thus,
\begin{align}\label{slope}
	\frac{\mathrm{d} u}{\mathrm{d} \rho}= - \frac{\gamma \widetilde{u}_0 p_0}{\rho^2} \left(\frac{1}{\rho} - \frac{1}{\rho^*}\right)^{\gamma-1} u^2 \left( \frac{1}{u} - \frac{1}{u^*}\right)^{2} < 0.
\end{align}
Consequently, the rarefaction and shock wave curves are both monotonically decreasing in the $(u, \rho)$ plane.

In the next part, we commence on the Riemann problem for \eqref{model1}. The Riemann solution plays an essential role in constructing the Riemann solvers to design numerical approaches, including the Godunov and Glimm schemes, as we will seen below. 

\subsection{The Riemann problem}
We study the Riemann problem associated with system \eqref{model1}, subject to the following initial conditions
\begin{align}
	(\rho, u) = 
	\left\{
	\begin{aligned}
		& (\rho_-, u_-), \quad x<0,\\
		& (\rho_+, u_+), \quad x>0,
	\end{aligned}
	\right.
\end{align}
in which $\rho_\pm$ and $ u_\pm$ are assigned constants.

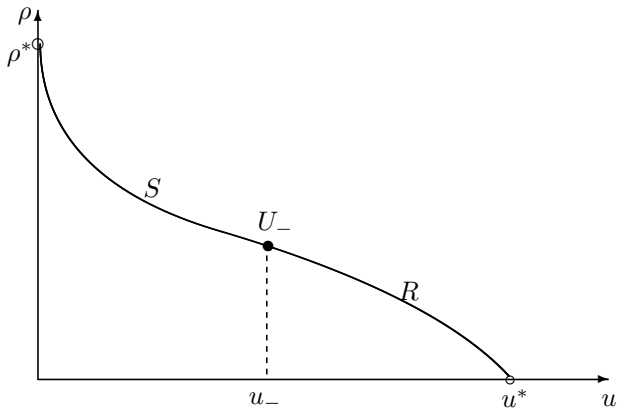
\begin{figure}[htbp]
	\centering
	\unitlength 0.75mm % = 2.845pt
	\linethickness{0.6pt}
	\ifx\plotpoint\undefined\newsavebox{\plotpoint}\fi % GNUPLOT compatibility
	\begin{picture}(100,70)(50,30)
		\put(54.5,37){\vector(1,0){100.75}}
		\put(54.562,36.909){\vector(0,1){64.878}}
		\qbezier(54.954,95.804)(55.116,72.056)(87.05,62.897)
		\qbezier(87.05,62.897)(124.334,51.874)(137.627,37.608)
		\put(54.555,95.88){\circle{1.734}}
		\put(137.875,36.875){\circle{1.458}}
		\put(95.156,60.465){\circle*{1.834}}
		%\dashline{1}(94.993,60.789)(94.993,36.96)
		\put(94.923,60.719){\line(0,-1){.9929}}
		\put(94.923,58.733){\line(0,-1){.9929}}
		\put(94.923,56.747){\line(0,-1){.9929}}
		\put(94.923,54.762){\line(0,-1){.9929}}
		\put(94.923,52.776){\line(0,-1){.9929}}
		\put(94.923,50.79){\line(0,-1){.9929}}
		\put(94.923,48.804){\line(0,-1){.9929}}
		\put(94.923,46.819){\line(0,-1){.9929}}
		\put(94.923,44.833){\line(0,-1){.9929}}
		\put(94.923,42.847){\line(0,-1){.9929}}
		\put(94.923,40.861){\line(0,-1){.9929}}
		\put(94.923,38.875){\line(0,-1){.9929}}
		%\end
		\put(155.297,33.394){\makebox(0,0)[cc]{$u$}}
		\put(138.762,33.718){\makebox(0,0)[cc]{$u^*$}}
		\put(94.669,33.232){\makebox(0,0)[cc]{$u_-$}}
		\put(96.452,64.194){\makebox(0,0)[cc]{$U_-$}}
		\put(51.387,94.021){\makebox(0,0)[cc]{$\rho^*$}}
		\put(52.198,100.829){\makebox(0,0)[cc]{$\rho$}}
		\put(74.73,70.516){\makebox(0,0)[cc]{$S$}}
		\put(119.958,52.522){\makebox(0,0)[cc]{$R$}}
	\end{picture}
	\vspace{-0.3cm}
	\caption{Shock and rarefaction wave curves}
	\label{rhou}
	\vspace{-0.2cm}
\end{figure}

From Eq.~\eqref{slope}, the curves $R(U, U_-)$ and $S(U, U_-)$ are depicted in Fig.~\ref{rhou}. Based on these curves, the problem can be classified into two distinct cases.

\begin{figure}[htbp]
	\centering
	\unitlength 0.65mm % = 2.845pt
	\linethickness{0.6pt}
	\ifx\plotpoint\undefined\newsavebox{\plotpoint}\fi % GNUPLOT compatibility
\begin{picture}(100,70)(100,40)
\put(70,50){\vector(1,0){70}}
\put(160,50){\vector(1,0){70}}
\put(105,50){\vector(0,1){50}}
\put(195,50){\vector(0,1){50}}
%\emline(104.961,50.027)(81.787,86.558)
\multiput(104.961,50.027)(-.0337308627,.05317571296){687}{\line(0,1){.05317571296}}
%\end
%\emline(104.824,49.89)(128.134,84.923)
\multiput(104.824,49.89)(.03373287271,.05069794319){691}{\line(0,1){.05069794319}}
%\end
%\emline(162.893,23.855)(163.166,23.991)
\multiput(162.893,23.855)(.054525,.027262){5}{\line(1,0){.054525}}
%\end
%\emline(194.849,49.928)(168.426,77.323)
\multiput(194.849,49.928)(-.03370266978,.03494325885){784}{\line(0,1){.03494325885}}
%\end
%\emline(195.011,49.928)(171.344,87.212)
\multiput(195.011,49.928)(-.03371386217,.05311087877){702}{\line(0,1){.05311087877}}
%\end
%\emline(194.686,50.252)(178.314,96.29)
\multiput(194.686,50.252)(-.0336882047,.0947272292){486}{\line(0,1){.0947272292}}
%\end
%\emline(195.011,50.09)(215.436,86.563)
\multiput(195.011,50.09)(.033704717,.0601869946){606}{\line(0,1){.0601869946}}
%\end
\put(139.733,46.848){\makebox(0,0)[cc]{$x$}}
\put(229.701,46.362){\makebox(0,0)[cc]{$x$}}
\put(102.449,98.883){\makebox(0,0)[cc]{$t$}}
\put(192.417,98.073){\makebox(0,0)[cc]{$t$}}
\put(80.403,88.888){\makebox(0,0)[cc]{$S$}}
\put(129.197,86.725){\makebox(0,0)[cc]{$J$}}
\put(85.267,60.465){\makebox(0,0)[cc]{$U_-$}}
\put(100.774,71.974){\makebox(0,0)[cc]{$U_1$}}
\put(123.037,59.816){\makebox(0,0)[cc]{$U_+$}}
\put(169.398,89.643){\makebox(0,0)[cc]{$R$}}
\put(217.219,91.102){\makebox(0,0)[cc]{$J$}}
\put(171.344,59.978){\makebox(0,0)[cc]{$U_-$}}
\put(213.977,60.303){\makebox(0,0)[cc]{$U_+$}}
\put(198.415,73.757){\makebox(0,0)[cc]{$U_2$}}
\put(105,47){\makebox(0,0)[cc]{$x=0$}}
\put(195,47){\makebox(0,0)[cc]{$x=0$}}
\end{picture}
	\vspace{-0.3cm}
	\caption{ The Riemann solution of case 1 and case 2.}
	\label{srj}
	\vspace{-0.2cm}
\end{figure}
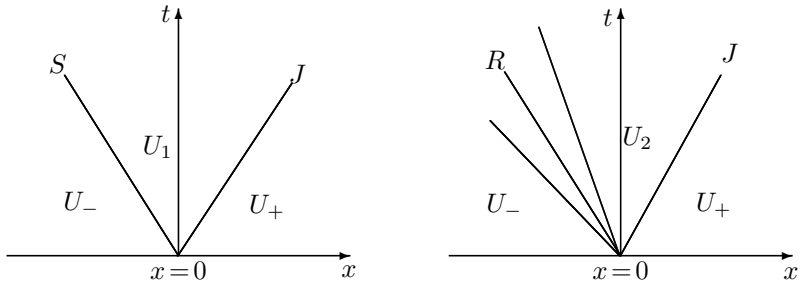

\paragraph{Case 1} $0<u_+<u_-$. Let $U_1$ denote the intersection point $J(U_+, U)\cap S(U, U_-) $. Then, the Riemann solution is given by (see Fig.~\ref{srj}, left)
$$S(U_1, U_-) \to J(U_+, U_1).$$ 

\paragraph{Case 2} $u_-<u_+<u^*$. Let $U_2$ denote the intersection point $ J(U_+, U)\cap R(U, U_-)$. In this case, the Riemann solution is expressed as (see Fig.~\ref{srj}, right)
$$R(U_2, U_-) \to J(U_+, U_2).$$ 

{\bf Remark.} When the ``vacuum" state occurs (i.e., $\rho \to 0$), the traffic velocity will approach its limit value $u^*$, which represents that in the totally free flow the vehicle will take the maximum allowed speed. Similarly, in the congestion flow, $\rho \to \rho^*$, thus the velocity $u\to 0$.

\section{Numerical approximation}\label{sec3}

The Godunov scheme for solving hyperbolic conservation laws is given by \cite{Godunov_1959}:
\begin{align}
	U_j^{n+1} = U_j^n - \frac{\Delta t}{\Delta x} \left( F_{j+1/2}^n - F_{j-1/2}^n \right), \quad j \in \mathbb{Z},
\end{align}
where $U_j^n$ denotes the cell-averaged conserved variables in the $j$-th cell at time step $n$, whereas the numerical flux $F_{j+1/2}^n$ is determined by solving the Riemann problem at the interface separating neighboring cells. Specifically,
\begin{align}
	F_{j+1/2}^n = F\left(U_r(0^-; U_j^n, U_{j+1}^n)\right),
\end{align}
where $U_r(0^-; U_j^n, U_{j+1}^n)$ corresponds to the left-hand limit at $x = x_{j+1/2}$, $t = 0$, derived from the local Riemann problem with initial data $U_j^n$ and $U_{j+1}^n$.

Although the Godunov scheme is robust in capturing shock and rarefaction waves, it can suffer from spurious oscillations near contact discontinuities when it is applied to the ARZ traffic flow model as pointed out in \cite{Christophe_2007_Godunov_Gimm}. To address this problem, the authors proposed a transport-equilibrium scheme that combines the Godunov scheme with the Glimm scheme that has infinite resolution in capturing the contact discontinuities. This hybrid approach significantly reduces numerical diffusion and suppresses nonphysical oscillations, thereby enhancing the fidelity and stability of numerical simulations in congested traffic scenarios.

To apply this strategy effectively, we follow the transport-equilibrium framework proposed in \cite{Christophe_2007_Godunov_Gimm}. Specifically, given the cell states $U_{j-1}^n$, $U_j^n$, and $U_{j+1}^n$, the update of $U_j^{n+1}$ within the spatial range $[x_{j-1}, x_{j+1}]$ is performed according to two successive steps described below.

we focus on the spatial interval $[x_{j-1}, x_{j+1}]$, the update of $U_j^{n+1}$ is carried out according to the following two steps. 

\begin{figure}[htbp]
	\centering
	\unitlength 0.5mm % = 2.845pt
	\linethickness{0.6pt}
	\ifx\plotpoint\undefined\newsavebox{\plotpoint}\fi % GNUPLOT compatibility
	\begin{picture}(100,50)(100,20)
		\put(70,30){\line(1,0){160}}
		\put(70,30){\line(0,1){2}}
		\put(110,30){\line(0,1){2}}
		\put(150,30){\line(0,1){2}}
		\put(190,30){\line(0,1){2}}
		\put(230,30){\line(0,1){2}}
		\put(70,25){\makebox(0,0)[cc]{$x_{j-1}$}}
		\put(150,25){\makebox(0,0)[cc]{$x_j$}}
		\put(190,25){\makebox(0,0)[cc]{$x_{j+1/2}$}}
		\put(230,25){\makebox(0,0)[cc]{$x_{j+1}$}}
		\put(110,25){\makebox(0,0)[cc]{$x_{j-1/2}$}}
		\put(70,64){\makebox(0,0)[cc]{$S/R$}}
		\put(145,64){\makebox(0,0)[cc]{$J$}}
		\put(160,67){\makebox(0,0)[cc]{$S/R$}}
		\put(222,67){\makebox(0,0)[cc]{$J$}}
		\put(75,40){\makebox(0,0)[cc]{$U_{j-1}^n$}}
		\put(150,40){\makebox(0,0)[cc]{$U_j^n$}}
		\put(215,40){\makebox(0,0)[cc]{$U_{j+1}^n$}}
		\put(110,58){\makebox(0,0)[cc]{$U^*(U_{j-1}^n, U_j^n)$}}
		\put(193,60){\makebox(0,0)[cc]{$U^*(U_j^n,U_{j+1}^n)$}}
		%\emline(109.956,29.876)(141.776,58.69)
		\multiput(109.956,29.876)(.03721650329,.03370161131){855}{\line(1,0){.03721650329}}
		%\end
		%\emline(190.39,29.876)(165.818,62.403)
		\multiput(190.39,29.876)(-.03370671645,.04461896279){729}{\line(0,1){.04461896279}}
		%\end
		%\emline(190.39,30.052)(156.802,61.519)
		\multiput(190.39,30.052)(-.03599988662,.03372620957){933}{\line(-1,0){.03599988662}}
		%\end
		%\emline(190.214,30.229)(219.028,62.756)
		\multiput(190.214,30.229)(.03370161131,.03804353669){855}{\line(0,1){.03804353669}}
		%\end
		%\emline(109.956,29.522)(80.434,58.16)
		\multiput(109.956,29.522)(-.03477266341,.0337315657){849}{\line(-1,0){.03477266341}}
		%\end
		%\emline(109.779,29.522)(71.595,57.983)
		\multiput(109.779,29.522)(-.045241863,.03372194418){844}{\line(-1,0){.045241863}}
		%\end
	\end{picture}
	\vspace{-0.3cm}
	\caption{The Riemann solution in $\left[x_{j-1}, x_{j+1}\right]$.}
	\label{step1}
	\vspace{-0.5cm}
\end{figure}
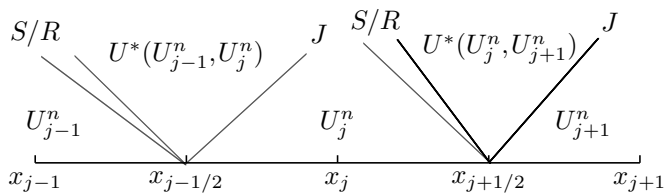

\paragraph{Step 1: Resolution of contact discontinuities ($t^n \to t^{n+1/2}$)}
In this step, special care is taken to accurately track the propagation of contact discontinuities. In the region $\left[x_{j-1}, x_{j+1}\right]$, the intermediate states $U^*(U_{j-1}^n, U_j^n)$ and $U^*(U_j^n, U_{j+1}^n)$ are obtained from solving the corresponding local Riemann problems. The contact discontinuities propagate rightward with positive speeds $u_j^n$ and $u_{j+1}^n$. The wave configuration, consisting of both first-family (shock/rarefaction) and second-family (contact) waves, is schematically depicted in Fig.~\ref{step1}.

Note that, although solving the full Riemann problem to obtain the intermediate states, the reconstruction procedure specifically utilizes only the contact discontinuities, without explicitly representing the rarefaction and shock waves. Then we construct a piecewise constant approximation, defined as follows:
\begin{align}
	U(x,t) = 
	\begin{cases}
		U_{j-1}^n, & \text{if } x \in \left[x_{j-1},\ x_{j - 1/2}\right), \\
		U^*(U_{j-1}^n, U_j^n), & \text{if } x \in \left[x_{j - 1/2},\ x_{j-1/2}+u_j^n(t-t^n)\right), \\
		U_j^n, & \text{if } x \in \left[x_{j-1/2}+u_j^n(t-t^n),\ x_{j+1/2}\right), \\
		U^*(U_j^n, U_{j+1}^n), & \text{if } x \in \left[x_{j+1/2},\ x_{j+1/2}+u_{j+1}^n(t-t^n)\right), \\
		U_{j+1}^n, & \text{if } x \in \left[x_{j+1/2}+u_{j+1}^n(t-t^n),\ x_{j+1}\right].
	\end{cases}
\end{align}

\begin{figure}[htbp]
	\centering
	\unitlength 0.7mm % = 2.845pt
	\linethickness{0.6pt}
	\ifx\plotpoint\undefined\newsavebox{\plotpoint}\fi % GNUPLOT compatibility
	\begin{picture}(100,35)(150,60)
		\put(140,70){\line(1,0){120}}
		\put(140,70){\line(0,1){2}}
		\put(200,70){\line(0,1){2}}
		\put(260,70){\line(0,1){2}}
		\put(200,65){\makebox(0,0)[cc]{$x_j$}}
		\put(140,65){\makebox(0,0)[cc]{$x_{j-1}$}}
		\put(170,65){\makebox(0,0)[cc]{$x_{j-1/2}$}}
		\put(230,65){\makebox(0,0)[cc]{$x_{j+1/2}$}}
		\put(260,65){\makebox(0,0)[cc]{$x_{j+1}$}}
		%\dashline{1}(167.338,68.113)(167.338,93.34)
		\put(170,69.983){\line(0,1){.9703}}
		\put(170,71.923){\line(0,1){.9703}}
		\put(170,73.864){\line(0,1){.9703}}
		\put(170,75.804){\line(0,1){.9703}}
		\put(170,77.745){\line(0,1){.9703}}
		\put(170,79.685){\line(0,1){.9703}}
		\put(170,81.626){\line(0,1){.9703}}
		\put(170,83.567){\line(0,1){.9703}}
		\put(170,85.507){\line(0,1){.9703}}
		\put(170,87.448){\line(0,1){.9703}}
		\put(170,89.388){\line(0,1){.9703}}
		\put(170,90){\line(0,1){.9703}}
		%\end
		%\dashline{1}(228.934,68.533)(228.934,91.868)
		\put(230,70){\line(0,1){.9723}}
		\put(230,72.352){\line(0,1){.9723}}
		\put(230,74.296){\line(0,1){.9723}}
		\put(230,76.241){\line(0,1){.9723}}
		\put(230,78.186){\line(0,1){.9723}}
		\put(230,80.13){\line(0,1){.9723}}
		\put(230,82.075){\line(0,1){.9723}}
		\put(230,84.019){\line(0,1){.9723}}
		\put(230,85.964){\line(0,1){.9723}}
		\put(230,87.908){\line(0,1){.9723}}
		\put(230,90){\line(0,1){.9723}}
		%\end
		\put(155,83){\makebox(0,0)[cc]{$U_{j-1}^n$}}
		\put(185,85){\makebox(0,0)[cc]{\footnotesize $U^*(U_{j-1}^n, U_j^n)$}}
		\put(207,76){\makebox(0,0)[cc]{$U_j^n$}}
		\put(247,83){\makebox(0,0)[cc]{$U_{j+1}^n$}}
		\put(212,88){\makebox(0,0)[cc]{$J$}}
		%\emline(167.339,68.323)(200.975,88.084)
		\multiput(170,70)(.0673994161,.026722157){586}{\line(1,0){.0573994161}}
		%\end
	\end{picture}
	\vspace{-0.3cm}
	\caption{The wave structure in $\left[x_{j-1}, x_{j+1}\right]$.}
	\label{step1.2}
	\vspace{-0.5cm}
\end{figure}
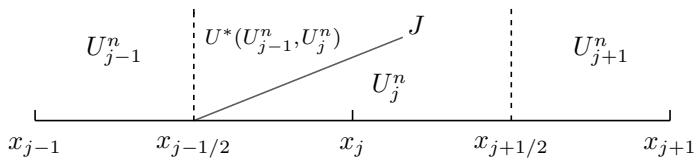

Du to the fact that the contact discontinuities propagate strictly to the right, these waves transport information into the target cell $C_j = [x_{j-1/2}, x_{j+1/2}]$ from its left neighbors, while remaining unaffected by information from the right. Based on this, the local wave structure within the region $\left[x_{j-1}, x_{j+1}\right]$, can be simplified as illustrated in Fig.~\ref{step1.2}.

To accurately resolve the contact discontinuities within the target cell $C_j$, Glimm’s random sampling strategy based on the van der Corput sequence $a_n$ is applied to enhance resolution. Specifically, the updated state is formulated as:
\begin{align} \label{half}
	U_j^{n+\frac{1}{2}} =  
	\begin{cases} 
		U^*(U_{j-1}^n, U_j^n),\  &\text{if} \quad a_{n+1} \in \left(0,\ \frac{\Delta t}{\Delta x} u_j^n\right), \\ 
		U_j^n,\  &\text{if} \quad a_{n+1} \in \left[ \frac{\Delta t}{\Delta x} u_j^n,\ 1 \right).
	\end{cases} 
\end{align} 
 %and $a_n = \sum_{k=0}^m i_k 2^{-(k+1)}$ represents the van der Corput low-discrepancy sequence. 

Accordingly, the piecewise constant approximation of the solution in the region $\left[x_{j-1}, x_{j+1}\right]$ is refined to incorporate the sampled intermediate state:
\begin{align}
	\widetilde{U}(x,t) = 
	\begin{cases}
		U_{j-1}^n, & \text{if } x \in \left[x_{j-1},\ x_{j - 1/2}\right), \\
		U_j^{n+\frac{1}{2}}, & \text{if } x \in \left[x_{j - 1/2}, x_{j+1/2}\right), \\
		U_{j+1}^n, & \text{if } x \in \left[x_{j+1/2},\ x_{j+1}\right].
	\end{cases}
\end{align}

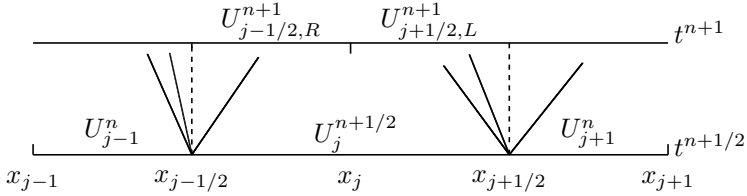
\begin{figure}[htbp]
	\centering
	\unitlength 0.7mm % = 2.845pt
	\linethickness{0.6pt}
	\ifx\plotpoint\undefined\newsavebox{\plotpoint}\fi % GNUPLOT compatibility
	\begin{picture}(100,40)(150,60)
		\put(140,70){\line(1,0){120}}
		\put(140,91){\line(1,0){120}}
		\put(140,70){\line(0,1){2}}
		%\put(200,70){\line(0,1){2}}
		\put(260,70){\line(0,1){2}}
		\put(200,89){\line(0,1){2}}
		\put(200,65){\makebox(0,0)[cc]{$x_j$}}
		\put(140,65){\makebox(0,0)[cc]{$x_{j-1}$}}
		\put(170,65){\makebox(0,0)[cc]{$x_{j-1/2}$}}
		\put(230,65){\makebox(0,0)[cc]{$x_{j+1/2}$}}
		\put(260,65){\makebox(0,0)[cc]{$x_{j+1}$}}
%\dashline{1}(167.338,68.113)(167.338,93.34)
\put(170,69.983){\line(0,1){.9703}}
\put(170,71.923){\line(0,1){.9703}}
\put(170,73.864){\line(0,1){.9703}}
\put(170,75.804){\line(0,1){.9703}}
\put(170,77.745){\line(0,1){.9703}}
\put(170,79.685){\line(0,1){.9703}}
\put(170,81.626){\line(0,1){.9703}}
\put(170,83.567){\line(0,1){.9703}}
\put(170,85.507){\line(0,1){.9703}}
\put(170,87.448){\line(0,1){.9703}}
\put(170,89.388){\line(0,1){.9703}}
\put(170,90){\line(0,1){.9703}}
%\end
%\dashline{1}(228.934,68.533)(228.934,91.868)
\put(230,70){\line(0,1){.9723}}
\put(230,72.352){\line(0,1){.9723}}
\put(230,74.296){\line(0,1){.9723}}
\put(230,76.241){\line(0,1){.9723}}
\put(230,78.186){\line(0,1){.9723}}
\put(230,80.13){\line(0,1){.9723}}
\put(230,82.075){\line(0,1){.9723}}
\put(230,84.019){\line(0,1){.9723}}
\put(230,85.964){\line(0,1){.9723}}
\put(230,87.908){\line(0,1){.9723}}
\put(230,90){\line(0,1){.9723}}
%\end
%x_j处的虚线
%\dashline{1}(228.934,68.533)(228.934,91.868)
%\put(200,70){\line(0,1){.9723}}
%\put(200,72.352){\line(0,1){.9723}}
%\put(200,74.296){\line(0,1){.9723}}
%\put(200,76.241){\line(0,1){.9723}}
%\put(200,78.186){\line(0,1){.9723}}
%\put(200,80.13){\line(0,1){.9723}}
%\put(200,82.075){\line(0,1){.9723}}
%\put(200,84.019){\line(0,1){.9723}}
%\put(200,85.964){\line(0,1){.9723}}
%\put(200,87.908){\line(0,1){.9723}}
%\put(200,90){\line(0,1){.9723}}
%\end
		\put(155,74){\makebox(0,0)[cc]{$U_{j-1}^n$}}
		\put(245,74){\makebox(0,0)[cc]{$U_{j+1}^n$}}
		\put(201,74){\makebox(0,0)[cc]{$U_j^{n+1/2}$}}
		\put(185,95){\makebox(0,0)[cc]{$U_{j-1/2,R}^{n+1}$}}
		\put(215,95){\makebox(0,0)[cc]{$U_{j+1/2,L}^{n+1}$}}
		\put(266,93){\makebox(0,0)[cc]{$t^{n+1}$}}
		\put(268,72){\makebox(0,0)[cc]{$t^{n+1/2}$}}
		%\emline(167.55,68.113)(163.345,87.244)
		\multiput(170,70)(-.033636058,.153044063){125}{\line(0,1){.153044063}}
		%\end
		%\emline(167.339,68.113)(158.93,87.033)
		\multiput(170,70)(-.0336360579,.0756811302){250}{\line(0,1){.0756811302}}
		%\end
		%\emline(167.55,68.533)(180.163,86.823)
		\multiput(170,70)(.0337259938,.0489026911){374}{\line(0,1){.0489026911}}
		%\end
		%\emline(228.725,68.533)(221.157,88.715)
		\multiput(230,70)(-.033636058,.083696154){225}{\line(0,1){.089696154}}
		%\end
		%\emline(228.725,68.744)(216.322,87.664)
		\multiput(230,70)(-.0337046096,.0454138113){368}{\line(0,1){.0514138113}}
		%\end
		%\emline(228.935,68.323)(242.81,85.562)
		\multiput(230,70)(.0336768783,.04184097){412}{\line(0,1){.04184097}}
		%\end
	\end{picture}
	\vspace{-0.3cm}
	\caption{Staggered grid update schematic.}
	\label{step2}
	\vspace{-0.5cm}
\end{figure}

\paragraph{Step 2: Resolution of Nonlinear Waves ($t^{n+1/2} \to t^{n+1}$)}
This step advances the solution by capturing nonlinear waves (shock and rarefaction waves). To compute the updated cell-average value $U_j^{n+1}$, we adopt the finite volume formulation, where the Godunov method yields the equivalence:
\begin{align}\label{u_j^n+1}
	U_j^{n+1} = \frac{1}{2} \left( U_{j-1/2,R}^{n+1} + U_{j+1/2,L}^{n+1} \right).
\end{align}
Here, $U_{j-1/2,R}^{n+1}$ represents the average solution on the interval $[x_{j-1/2}, x_j]$, while $U_{j+1/2,L}^{n+1}$ corresponds to the average over $[x_j,, x_{j+1/2}]$. These values as shown in Fig.~\ref{step2}, and they incorporate the effect of wave propagating into the target cell $C_j$ from both the right and the left interfaces.

At the right interface $x_{j+1/2}$, only waves propagating to the left can enter cell $C_j$. Therefore, we integrate the Riemann solution corresponding to the states $U_j^{n+1/2}$ and $U_{j+1}^n$ within the domain $(x_j, x_{j+1/2})$, leading to
\begin{align}\label{right}
	\begin{aligned}
		U_{j+1/2,L}^{n+1} &= \frac{2}{\Delta x} \int_{x_j}^{x_{j+1/2}} U_r\left(\frac{x - x_{j+1/2}}{\Delta t}; U_j^{n+1/2}, U_{j+1}^n\right) \,dx \\
		&= U_j^{n+1/2} - \frac{2\Delta t}{\Delta x} \left( F(U_r(0^-; U_j^{n+1/2}, U_{j+1}^n)) - F(U_j^{n+1/2}) \right),
	\end{aligned}
\end{align}
where $U_r(\cdot\,;\,U_j^{n+1/2}, U_{j+1}^n)$ denotes the self-similar Riemann solution with the intermediate state at $x = x_{j+1/2}$ given by $U_r(0^-)$, which accounts for all left-going waves.

The treatment of the left interface $x_{j-1/2}$ demands special attention, as its dynamics depend critically on whether the 1-wave propagates into the target cell $C_j$ or not. The solution depends on the random sampling result, which leads to two main cases:

\paragraph{Case 1} If $U_j^{n+1/2} = U_j^n$,  this sampling result implies that the contact discontinuity has not yet propagated into the target cell $C_j$. However, the exact Riemann solution $U_r(\cdot; U_{j-1}^n, U_j^n)$ may still contain a contact discontinuity, giving rise to two possible subcases.

\paragraph{Subcase 1} No contact discontinuity exists in the Riemann solution, yet a 1-wave may still cross into the cell $C_j$. In this situation, the update state must take account for the contribution from the 1-wave. Since $U_j^{n+1/2} = U_j^n$, the intercell numerical flux can be computed alternatively by solving the equivalent Riemann problem $U_r(\cdot; U_{j-1}^n, U_j^{n+1/2})$.

\paragraph{Subcase 2} A contact discontinuity is present, which leads to $ U^*(U_{j-1}^n, U_j^n) \neq U_j^n $, with $U^*$ denoting the intermediate Riemann state. However, since the sampling step indicates that the discontinuity has not yet entered the cell $C_j$. Thus, the physical state within the target cell is governed by the right-going characteristics emanating from the contact discontinuity, i.e., it is still given by $U_j^n = U_j^{n+1/2}$.

\paragraph{Case 2.} If $U_j^{n+1/2} = U^*(U_{j-1}^n, U_j^n) \neq U_j^n$, this indicates that the contact discontinuity has propagated into the cell $C_j$. In this case, to properly account for the incoming 1-wave across the interface, the intercell numerical flux must be computed by solving the Riemann problem  $U_r(\cdot; U_{j-1}^n, U_j^{n+1/2})$.

From  subcase 2 of case 1, since $U_j^{n+1/2} = U_j^n$, the condition $U^*(U_{j-1}^n, U_j^n) \neq U_j^n$ is equivalent to $U^*(U_{j-1}^n, U_j^{n+1/2}) \neq U_j^{n+1/2}$. Both case 1 (subcase 1) and case 2 ultimately rely on solving the same Riemann problem $U_r(\cdot; U_{j-1}^n, U_j^{n+1/2})$, which enables us to unify the update procedure into the following two cases.

If $U^*(U_{j-1}^n, U_j^{n+1/2}) \neq U_j^{n+1/2}$, then no wave has entered the cell from the left. Therefore, the state can be directly taken as the cell value:
\begin{align}\label{left1}
	U_{j-1/2,R}^{n+1} = U_j^{n+1/2}.
\end{align}
Otherwise, the 1-wave has crossed into the cell. Analogous to the derivation at the right interface, we integrate the Riemann solution over the interval $(x_{j-1/2}, x_j)$, and have
\begin{align}\label{left2}
	U_{j-1/2,R}^{n+1} = U_j^{n+1/2} - \frac{2\Delta t}{\Delta x} \left( F(U_j^{n+1/2}) - F(U_r(0^-; U_{j-1}^n, U_j^{n+1/2})) \right).
\end{align}

Thus, by \eqref{u_j^n+1}, combining \eqref{right}, \eqref{left1} and \eqref{left2}, we have the following update solution
\begin{align} 
	U_j^{n+1} = \frac{1}{2} \left( U_{j,L}^{n+1} + U_{j,R}^{n+1}\right) = U_j^{n+\frac{1}{2}} - \frac{\Delta t}{\Delta x} \left( F_{j+1/2}^{n+1/2, L} - F_{j-1/2}^{n+1/2,R} \right),
\end{align}
where the numerical flux are defined according to
\begin{align*}
	&F_{j+1/2}^{n+1/2, L} = F\big(U_r(0^-; U_j^{n+1/2}, U_{j+1}^n)\big),\\
	&F_{j-1/2}^{n+1/2, R} = 
	\begin{cases}
		F(U_j^{n+1/2}), & \text{if } U^*(U_{j-1}^n, U_j^{n+1/2}) \neq U_j^{n+1/2},\\
		F\big(U_r(0^-; U_{j-1}^n, U_j^{n+1/2})\big), & \text{otherwise}.
	\end{cases}
\end{align*}

Now, we test some 1-D numerical examples to compare the present the model with other models, and also to see how the numerical method performs.  

{\em Test 1.} Firstly, we compare the three models: ARZ, MAR, and RARZ,  under the same initial conditions, which is given as
\begin{align}
	U(x,0) = \left\{
	\begin{aligned}
		& U_-=(\rho_-, u_-)=(0.4, 20),\quad x<1\\
		& U_+=(\rho_+, u_+)=(0.8, 16),\quad x>1.
	\end{aligned}
	\right.
\end{align}
\begin{figure}[h!]
	\centering
	\vspace{-0.5cm}
	\includegraphics[width=0.8\linewidth]{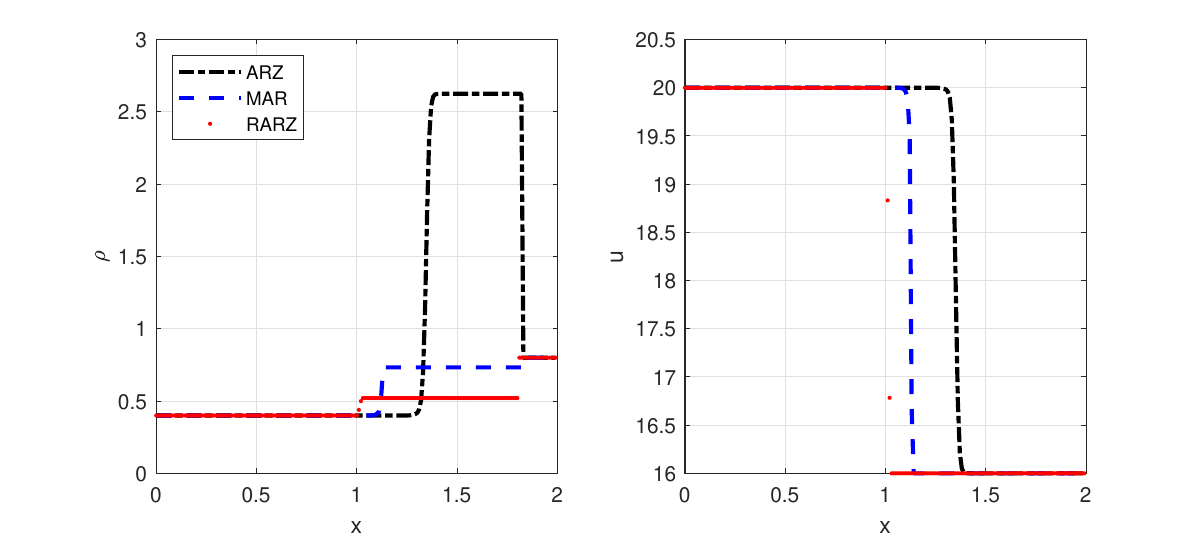}
	\caption{Test 1: the solution for density and velocity at $ t = 0.05\text{s} $.}
	\label{modelc}
	\vspace{-0.7cm}
\end{figure}

From Fig. 4.4, compared to the ARZ model, the MAR model has a ``quicker" response to the potential traffic jam. This is because the modified pressure has a higher sensitivity to the change of the traffic density.  Moreover, the present RARZ model  has a ``quicker" response than the MAR model, which demonstrates the present model can alleviate the congestion probability more efficiently. 
 
Now, we present numerical simulations of the Riemann solution for system \eqref{model1}. For the following numerical tests, the simulation is defined on the domain $[0,2]$ and features an initial discontinuity at the midpoint. We present a comparative analysis of the classical Godunov with  the hybrid Godunov–Glimm scheme. 

It is evident that the classical Godunov method introduces noticeable numerical diffusion, in contrast, the Godunov–Glimm scheme closely matches the exact solution. This hybrid method effectively preserves sharp discontinuities and captures wave interactions with high fidelity, demonstrating its ability in resolving fine-scale structures in traffic flow dynamics.

{\em Test 2.} The initial date is given as
\begin{align}
	U(x,0) =\left\{
	\begin{aligned}
		& U_-=(\rho_-, u_-)=(0.8, 22),\quad x<1\\
		& U_+=(\rho_+, u_+)=(0.6, 15),\quad x>1.
	\end{aligned}
	\right.
\end{align}
\begin{figure}[h!]
	\centering
	\vspace{-0.7cm}
	\includegraphics[width=0.8\linewidth]{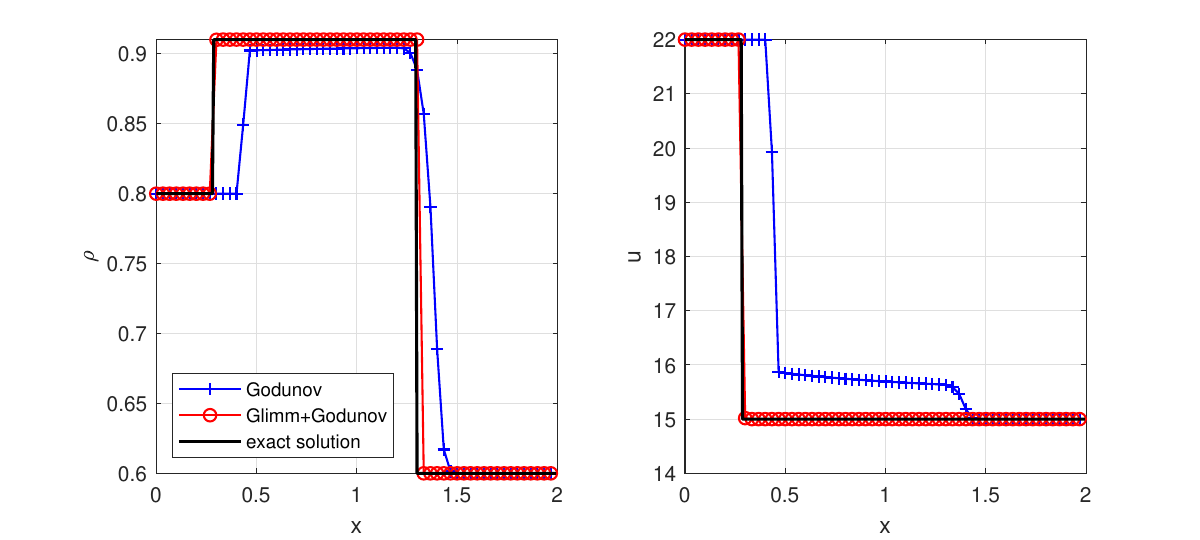}
	\caption{Test 2: the solution for density and velocity at $ t = 0.02\text{s} $.}
	\label{methodc1}
	\vspace{-0.7cm}
\end{figure}

The output time is taken at $ t = 0.02\text{s} $, the Riemann solution consists of a shock wave $S$, followed by a contact discontinuity $J$, see Fig.~\ref{methodc1}. We can see that the hybrid method captures the exact solution quite accurately than the Godunov scheme for both waves.

{\em Test 3.} The initial date is given as
\begin{align}
	U(x,0) =\left\{
	\begin{aligned}
		& U_-=(\rho_-, u_-)=(0.8, 16),\quad x<1\\
		& U_+=(\rho_+, u_+)=(0.6, 18),\quad x>1.
	\end{aligned}
	\right.
\end{align}

\begin{figure}[h!]
	\centering
	\vspace{-0.5cm}
	\includegraphics[width=0.8\linewidth]{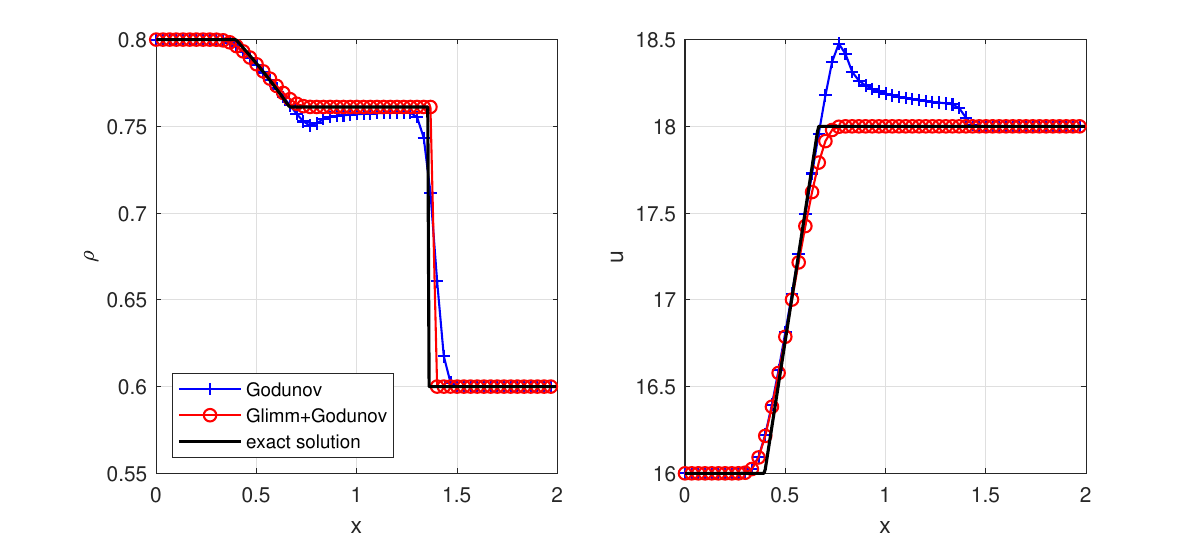}
	\caption{Test 3: the solution for density and velocity at $ t = 0.02\text{s} $.}
	\label{methodc2}
	\vspace{-0.1cm}
\end{figure}

The output time is taken at $ t = 0.02\text{s} $, the Riemann solution consists of a rarefaction wave $R$, followed by a contact discontinuity $J$, see Fig.~\ref{methodc2}. Again, we see that the hybrid method is in excellent agreement with the exact solution. 

{\em Test 4.} The initial date is given as
\begin{align}
	U(x,0) =\left\{
	\begin{aligned}
		& U_-=(\rho_-, u_-)=(0.8, 15),\quad x<1\\
		& U_+=(\rho_+, u_+)=(0.7, 15),\quad x>1.
	\end{aligned}
	\right.
\end{align}

\begin{figure}[h!]
	\centering
	\vspace{-0.1cm}
	\includegraphics[width=0.8\linewidth]{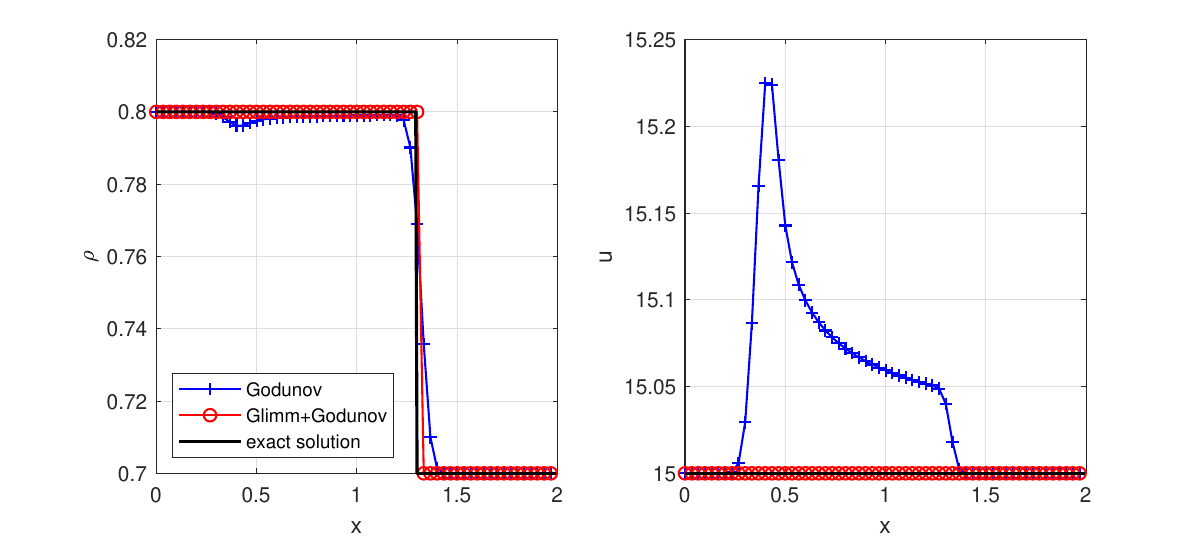}
	\caption{Test 4: the solution for density and velocity at $ t = 0.02\text{s} $.}
	\label{methodc3}
	\vspace{-0.9cm}
\end{figure}

The output time is taken at $ t = 0.02\text{s} $, the numerical result indicates that the flow pattern is characterized solely by a contact discontinuity $J$, which corresponds to the analytical Riemann solution.

\section{Two-dimensional RARZ traffic flow model}\label{sec4}
The present section formulates a 2D traffic flow framework as an extension of the 1D case, building on the framework introduced in \cite{Herty_2018_2D} when lane changing is considered. We derive a two-dimensional microscopic model by systematically incorporating fundamental traffic flow parameters, namely, $\rho^*$, $u^*$, and $v^*$ are defined as the maximum values of density, x-velocity, and y-velocity, respectively.  

Generalizing the one-dimensional definition to two dimensions, the local density centered at vehicle $i$ is given by:
$$\rho_i = \frac{\Delta X \Delta Y}{(x_{j(i)} - x_i) |y_{j(i)} - y_i|},$$
here, $j(i)$ denotes the one immediately ahead of vehicle $i$, $\Delta X$, $\Delta Y$ represent the longitudinal and lateral dimensions of the vehicle. 

The corresponding specific volume is defined as:
$$\tau _i = \frac{(x_{j(i)} - x_i) |y_{j(i)} - y_i| - d_x d_y }{\Delta X \Delta Y} =\frac{1}{\rho_i} - \frac{1}{\rho^*},$$
where $d_x, d_y$ represent the minimal admissible distance in the longitudinal and lateral directions and satisfy $\frac{d_x d_y}{\Delta X \Delta Y} = \frac{1}{\rho^*} $. 

Differentiating the specific volume and applying the density definition leads to:
\begin{align*}
	\dot{\tau}_i = \frac{u_{j(i)} - u_i}{\Delta X} + \frac{(v_{j(i)} - v_i) \left| y_{j(i)} - y_i\right|}{\Delta Y (y_{j(i)} - y_i)}.
\end{align*}

To characterize the desired motion states along the longitudinal ($x$) and lateral ($y$) directions, two auxiliary variables are introduced:
\begin{align*}
	&w_i = \widetilde{u}_i\, p(\tau_i) = \left(\frac{1}{u_i} - \frac{1}{u^*}\right)^{-1} \left(\frac{1}{\rho_i} - \frac{1}{\rho^*}\right)^{-\gamma},\\
	&\sigma_i = \widetilde{v}_i\, p(\tau_i) = \left(\frac{1}{v_i} - \frac{1}{v^*}\right)^{-1} \left(\frac{1}{\rho_i} - \frac{1}{\rho^*}\right)^{-\gamma}.
\end{align*}

Under the fundamental assumption that these desired velocities remain unchanged over time, we impose the conservation conditions
$$\dot{w}_i = 0, \qquad \dot{\sigma}_i = 0.$$

As a result of this construction, the dynamics can be represented by a two-dimensional microscopic framework of the follow-the-leader type, formulated as follows
\begin{align}\label{2Dm1}
	\left\{
	\begin{aligned}
		& \dot{x}_i = u_i,\\
		& \dot{y}_i = v_i,\\
		& \dot{u}_i = \frac{\gamma u_i (u^* - u_i)}{u^*} \left( \frac{u_{j(i)} - u_i}{x_{j(i)} - x_i - \frac{dl}{|y_{j(i)}-y_i|}} + \frac{v_{j(i)} - v_i}{y_{j(i)} - y_i - \frac{dl(y_{j(i)}-y_i)}{(x_{j(i)} - x_i)|y_{j(i)}-y_i|}} \right), \\
		& \dot{v}_i = \frac{\gamma v_i (v^* - v_i)}{v^*} \left( \frac{u_{j(i)} - u_i}{x_{j(i)} - x_i - \frac{dl}{|y_{j(i)}-y_i|}} + \frac{v_{j(i)} - v_i}{y_{j(i)} - y_i - \frac{dl(y_{j(i)}-y_i)}{(x_{j(i)} - x_i)|y_{j(i)}-y_i|}} \right).
	\end{aligned}
	\right.
\end{align}

It is worth noting that when $v_i = 0$, and assuming that
$$\lim\limits_{\substack{d_y \to 0^+ \\ y_{j(i)} - y_i \to 0^+ }} \frac{d_y}{|y_{j(i)}-y_i|} = 1,$$
 then the two-dimensional model naturally reduces to its one-dimensional counterpart, as presented in Eq.\eqref{mar}.

For analytical convenience, the two-dimensional system \eqref{2Dm1} admits the following equivalent representation
\begin{align}\label{2Dm2}
	\left\{
	\begin{aligned}
		& \dot{\tau}_i = \frac{u_{j(i)} - u_i}{\Delta X} + \frac{(v_{j(i)} - v_i) \left| y_{j(i)} - y_i \right|}{\Delta Y (y_{j(i)} - y_i)}, \\
		& \dot{w}_i = 0, \\
		& \dot{\sigma}_i = 0.
	\end{aligned}
	\right.
\end{align}

By taking the continuum limit, we obtain the corresponding macroscopic system in Lagrangian coordinates:
\begin{align}
	\partial _t \tau^L = \partial _X u^L + \partial _Y v^L, \quad \partial _t w^L =0, \quad \partial_t \sigma^L =0.
\end{align}
Finally, the transition from Lagrangian to Eulerian coordinates yields the following macroscopic equations:
\begin{align} \label{model2}
	\left\{
	\begin{aligned}
		& \rho_t + (\rho u)_x + (\rho v)_y = 0, \\
		& (\rho \widetilde{u} p)_t + (\rho u \widetilde{u} p)_x + (\rho v \widetilde{u} p)_y = 0, \\
		& (\rho \widetilde{v} p)_t + (\rho u \widetilde{v} p)_x + (\rho v \widetilde{v} p)_y = 0,
	\end{aligned}
	\right.
\end{align}
where $\widetilde{u} = \left(\frac{1}{u} - \frac{1}{u^*}\right)^{-1}, \widetilde{v} = \left(\frac{1}{v} - \frac{1}{v^*}\right)^{-1}, p=\left(\frac1\rho-\frac{1}{\rho^*}\right)^{-\gamma}$.

Denote the conservative variable vector by  $U = \left( \rho, \ \rho \widetilde{u} p, \ \rho \widetilde{v} p \right)^T$, the original system \eqref{model2} can be reformulated in conservative form:
$$ \partial_t U + \partial_x F(U) + \partial_y G(U) =0, $$
here, $F(U)$ and $G(U)$ denote the flux functions along the longitudinal $(x)$ and lateral $(y)$ directions, respectively.

We adopt the Strang splitting method which decomposes the two-dimensional problem into two 1-D systems at each direction.
We first consider the subsystem along the longitudinal $(x)$ direction
\begin{align}\label{x-direction}
	\partial_t U + \partial_x F(U) =0.
\end{align}
For this subsystem, the Jacobian matrix $\partial F / \partial U$ possesses three eigenvalues:
$$ \lambda_1= \lambda_2 = u , \quad \lambda_3=u - \frac{\gamma \rho^* u \left(u^* - u\right)}{u^*\left(\rho^* -\rho \right)}.$$
The associated right eigenvectors are
$$ \vec{r}_1 = \vec{r}_2 = (1, \ 0, \ 1)^T, \quad \vec{r}_3=\left(-\rho, \ \frac{\gamma \rho^* u \left(u^* - u\right)}{u^*\left(\rho^* -\rho \right)}, \ \frac{\gamma \rho^* v \left(v^* - v\right)}{v^*\left(\rho^* -\rho \right)}\right)^T.$$
The Riemann invariant $w_i$, corresponding to the $i$-th characteristic field $\lambda_i$, is
$$w_1 = u, \quad w_2 = \widetilde{u} p, \quad w_3 = \widetilde{v} p.$$
For system \eqref{x-direction}, the Rankine–Hugoniot conditions take the form
\begin{align*}
	\left\{
	\begin{aligned}
		& \sigma[\rho]=[\rho u] ,\\
		& \sigma[\rho \widetilde{u} p]=[\rho u \widetilde{u} p],\\
		& \sigma[\rho \widetilde{v} p]=[\rho u \widetilde{v} p].
	\end{aligned}
	\right.
\end{align*}
where $[f]:=f_1-f_0$, the conditions are given by $u = u_0$, $\widetilde{u}p = \widetilde{u}_0p_0$, and $\widetilde{v}p = \widetilde{v}_0p_0$.

The system in the $y$-direction is
$$\partial_t U + \partial_y G(U) = 0,$$
the same procedure can be applied to this system by transforming $x\to y$ and $\widetilde{u} \to \widetilde{v}$.

For each cell interface, the intermediate state $(\rho_M, u_M, v_M)$ is computed by enforcing consistency through Riemann invariants and the Rankine-Hugoniot conditions within a nonlinear framework, as follows
\begin{align*}
	\left\{
	\begin{aligned}
		&u_M = u_R,\\
		&\widetilde{u}_L p_L = \widetilde{u}_M p_M,\\
		&\widetilde{v}_L p_L = \widetilde{v}_M p_M.
	\end{aligned}
	\right.
\end{align*}

Subsequently, the full time step $\Delta t$ is implemented through three fractional stages, following the Strang splitting strategy. The numerical updates are carried out as follows
\begin{align*}
	& U_{i,j}^{n+1/3} = U_{i,j}^n - \frac{\Delta t}{2 \Delta x} \left( F(U_{i+1/2,j}^n) - F(U_{i-1/2,j}^n) \right), \\
	& U_{i,j}^{n+2/3} = U_{i,j}^{n+1/3} - \frac{\Delta t}{\Delta y} \left( G(U_{i,j+1/2}^{n+1/3}) - G(U_{i,j-1/2}^{n+1/3}) \right), \\
	& U_{i,j}^{n+1} = U_{i,j}^{n+2/3} - \frac{\Delta t}{2 \Delta x} \left( F(U_{i+1/2,j}^{n+2/3}) - F(U_{i-1/2,j}^{n+2/3}) \right),
\end{align*}
here, he cell-interface fluxes $F_{i\pm\frac{1}{2},j}$ and $G_{i,j\pm\frac{1}{2}}$ are computed via the HLL (Harten--Lax--van Leer) approximate Riemann solver. This fractional-step method achieves dimensional decoupling while maintaining conservation properties.

To validate the proposed model, we conduct two numerical tests by solving the Riemann problem for system \eqref{model2} on the computational domain $\Omega = [0,2] \times [0,2]$. Both tests implement piecewise-constant initial conditions with discontinuities along $x = 2$ and $y = 2$. A comparative view of Test 1, Test 2, and Test 3 is provided in Figs. \ref{fourS}, \ref{fourJ}, and \ref{2R2J}. Each composite figure combines a two-dimensional density contour plot (left) illustrating wave patterns, and a corresponding three-dimensional density surface plot (right) that reveals the amplitude variations.

{\em Test 1.} The initial date is given as
\begin{align*}
	(\rho_0, u_0 ,v_0)=
	\left\{
	\begin{aligned}
		&(\rho_1, u_1, v_1) = (0.4275, 0.5, 0.2), \quad (x,y)\in [1,2]\times [1,2]\\
		&(\rho_2, u_2, v_2) = (0.3199, 0.8, 0.2), \quad (x,y)\in [0,1]\times [1,2]\\
		&(\rho_3, u_3, v_3) = (0.2152, 0.8, 0.4), \quad (x,y)\in [0,1]\times [0,1]\\
		&(\rho_4, u_4, v_4) = (0.3033, 0.5, 0.4), \quad (x,y)\in [1,2]\times [0,1].
	\end{aligned}
	\right.
\end{align*}

\begin{figure}[h!]
	\centering
	\vspace{-0.5cm}
	\includegraphics[width=0.8\linewidth]{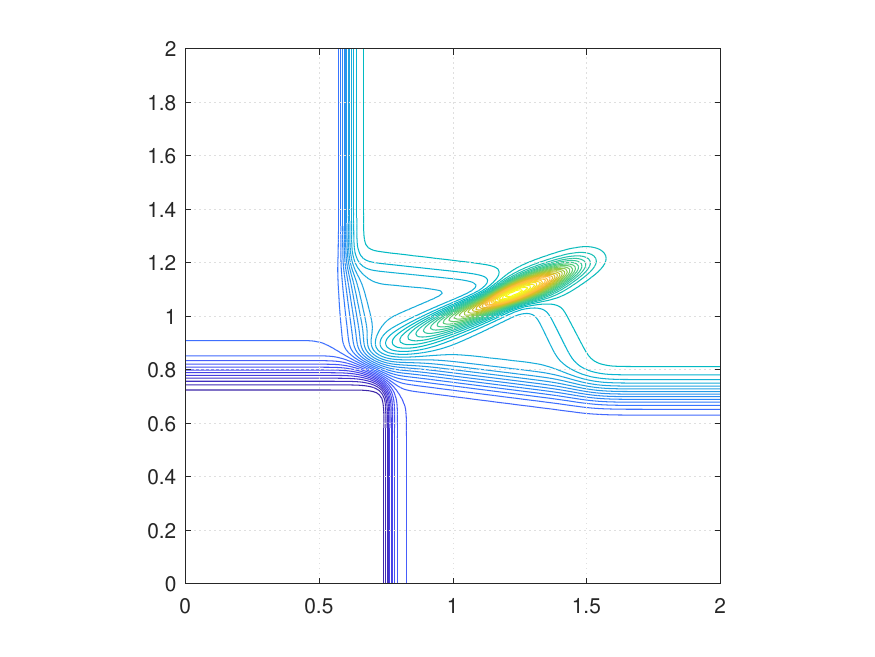}
		\caption{Density contours (HLL scheme)}
			\label{hll}
	\vspace{-0.1cm}
\end{figure}

\begin{figure}[h!]
	\centering
	\vspace{-0.5cm}
		\includegraphics[width=\linewidth]{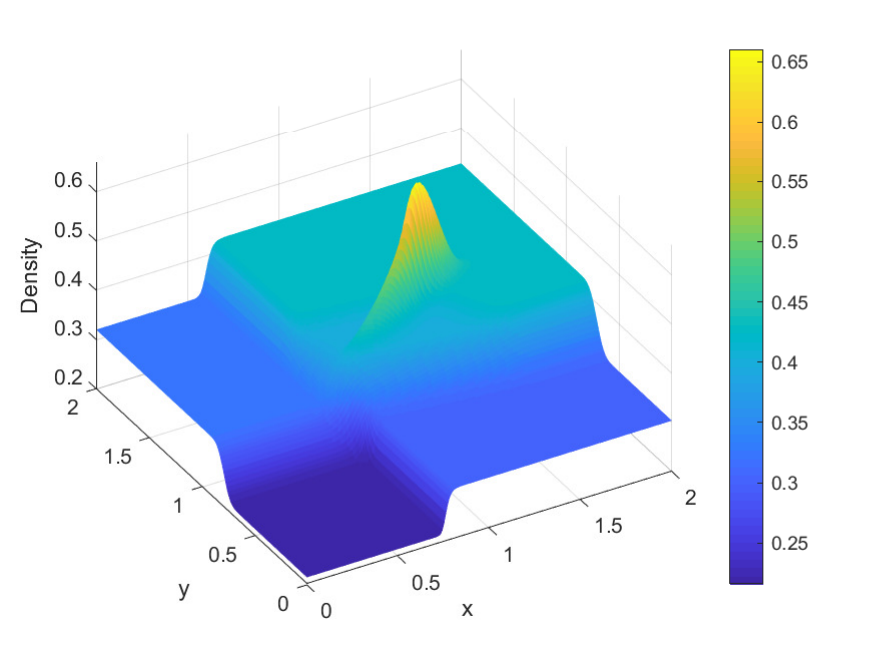}
		\caption{3D density field (HLL scheme)}
		\label{hll_3d}
	\vspace{-0.1cm}
\end{figure}

\begin{figure}[h!]
	\centering
	\vspace{-0.5cm}
		\includegraphics[width=\linewidth]{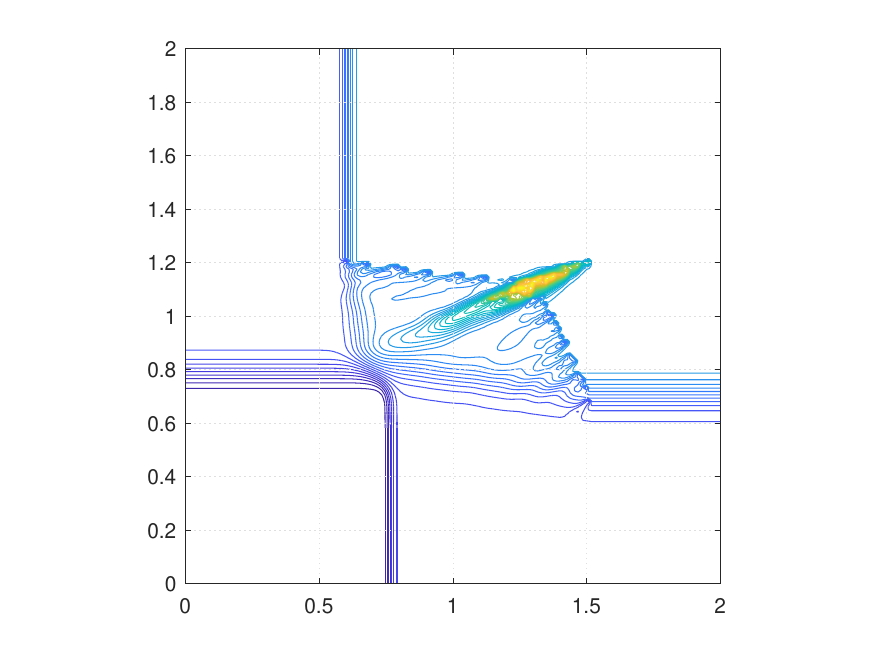}
		\caption{Density contours (Godunov-Glimm scheme)}
		\label{gg}
	\vspace{-0.1cm}
\end{figure}

\begin{figure}[h!]
	\centering
	\vspace{-0.5cm}
	\includegraphics[width=\linewidth]{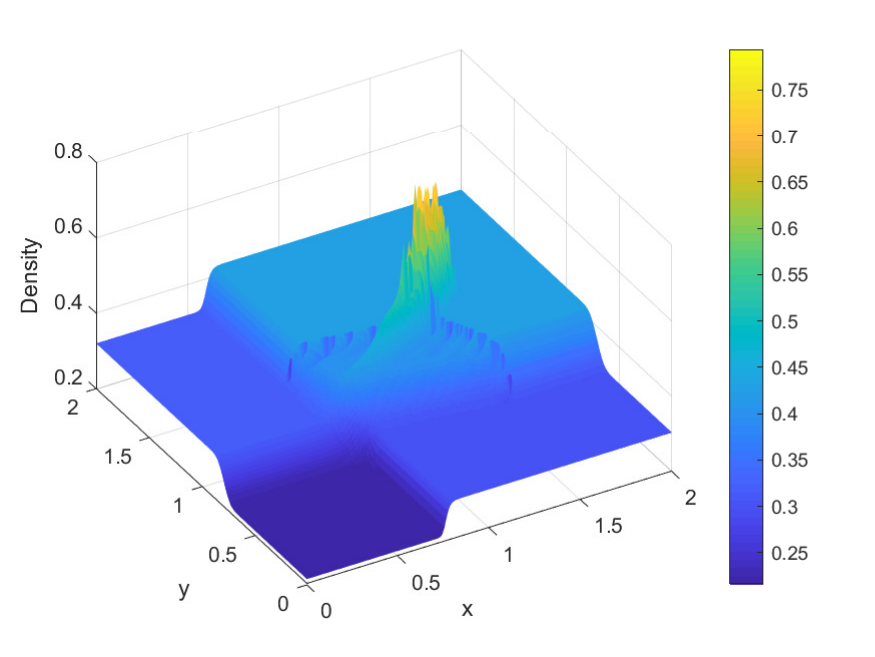}
		\caption{3D density field (Godunov-Glimm scheme)}
		\label{gg_3d}
	\vspace{-0.1cm}
\end{figure}

The numerical solution is output at $t = 1\,\text{s}$, where the flow field exhibits four shock waves as shown in Fig.~\ref{fourS}. Figs.~\ref{hll} and \ref{gg} present the shock wave contour plots obtained using the HLL scheme and the Godunov-Glimm hybrid method, respectively, while Figs.~\ref{hll_3d} and ~\ref{gg_3d} display their corresponding three-dimensional density field visualizations. Comparative analysis demonstrates that the Godunov-Glimm method achieves a good resolution, particularly in maintaining sharper contact discontinuities and reducing numerical dissipation.

{\em Test 2.} The initial date is given as
\begin{align*}
	(\rho_0, u_0 ,v_0)=
	\left\{
	\begin{aligned}
		&(\rho_1, u_1, v_1) = (0.2, 0.8, 0.3), \quad (x,y)\in [1,2]\times [1,2]\\
		&(\rho_2, u_2, v_2) = (0.3, 0.8, 0.5), \quad (x,y)\in [0,1]\times [1,2]\\
		&(\rho_3, u_3, v_3) = (0.4, 0.5, 0.5), \quad (x,y)\in [0,1]\times [0,1]\\
		&(\rho_4, u_4, v_4) = (0.5, 0.5, 0.3), \quad (x,y)\in [1,2]\times [0,1].
	\end{aligned}
	\right.
\end{align*}

\begin{figure}[h!]
	\centering
	\vspace{-0.5cm}
		\includegraphics[width=\linewidth]{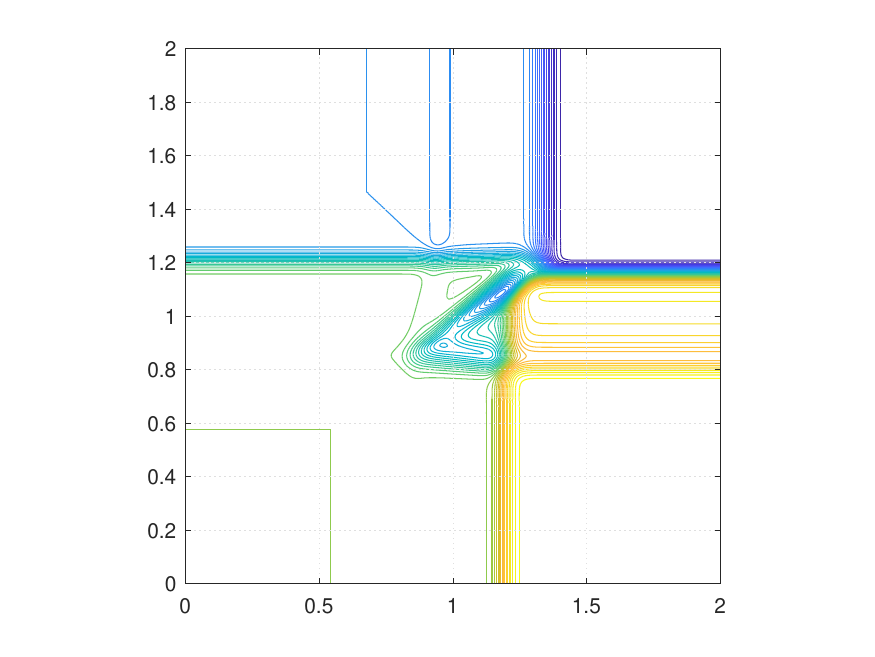}
		\caption{Density contours (HLL scheme)}
		\label{J_hll}
	\vspace{-0.1cm}
\end{figure}

\begin{figure}[h!]
	\centering
	\vspace{-0.5cm}
		\includegraphics[width=\linewidth]{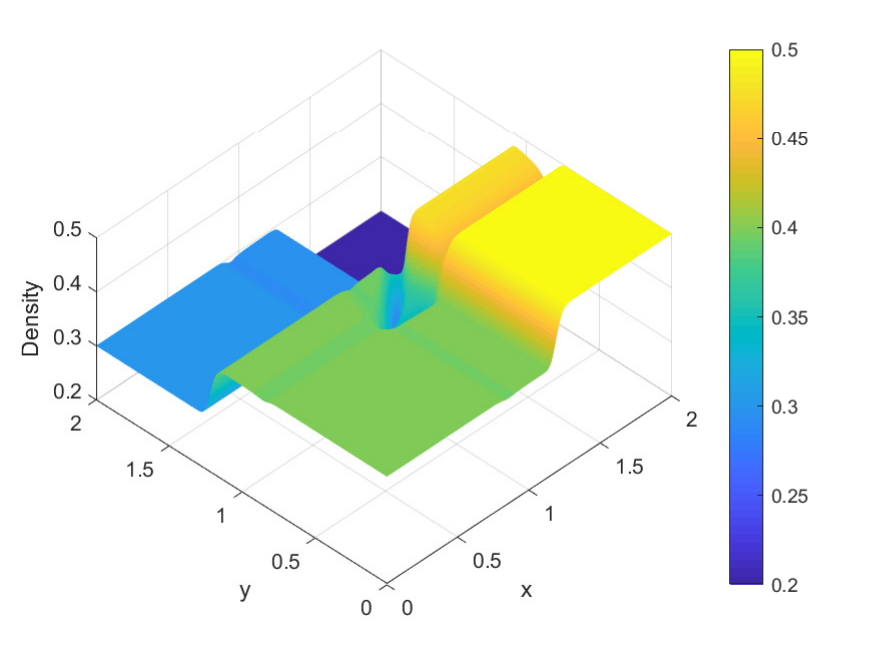}
		\caption{3D density field (HLL scheme)}
		\label{J_hll_3d}
	\vspace{-0.1cm}
\end{figure}

	\begin{figure}[h!]
	\centering
	\vspace{-0.5cm}
		\includegraphics[width=\linewidth]{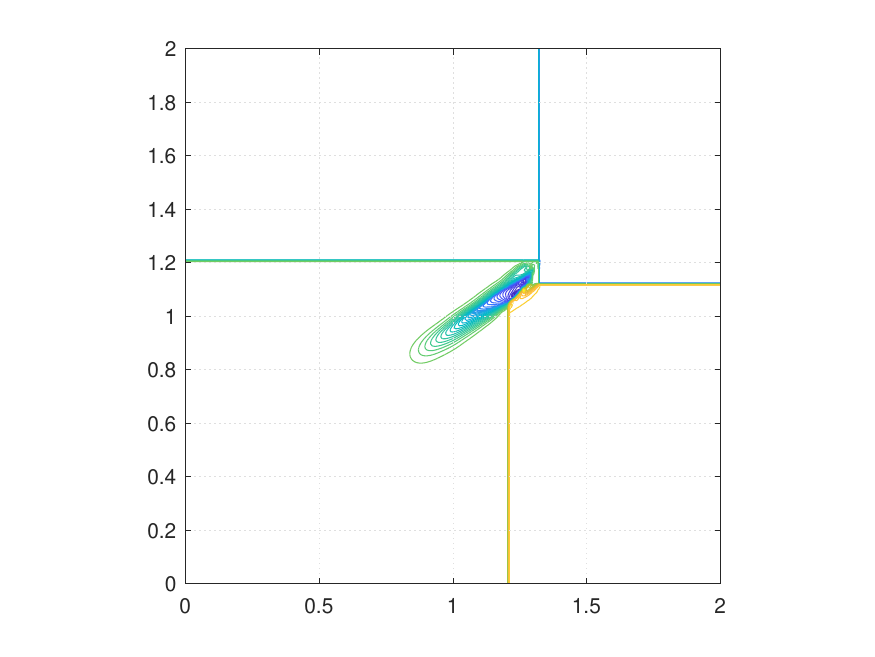}
		\caption{Density contours (Godunov-Glimm scheme)}
		\label{J_gg}
	\vspace{-0.1cm}
\end{figure}

\begin{figure}[h!]
	\centering
	\vspace{-0.5cm}
		\includegraphics[width=\linewidth]{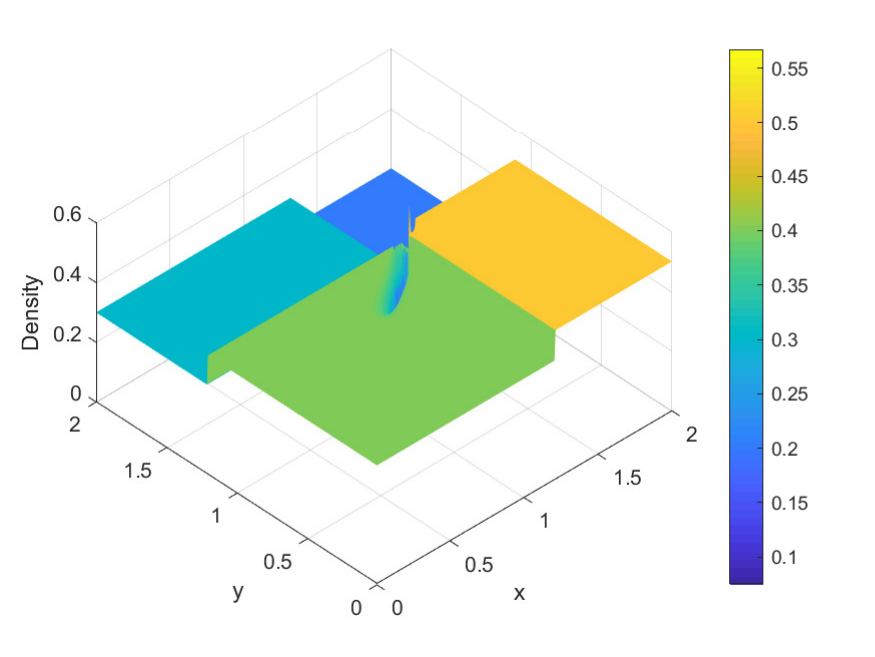}
		\caption{3D density field (Godunov-Glimm scheme)}
		\label{J_gg_3d}
	\vspace{-0.1cm}
\end{figure}

The numerical solution is output at the final simulation time $t = 0.4\,\text{s}$, where the flow field exhibits four contact discontinuities as shown in Fig.~\ref{fourJ}. One can see clearly that the result obtained by the HLL scheme is spurious, and the positions of the contact discontinuities are also not correct. While the Godunov-Glimm method can track the right positions of the contact discontinuities, the interaction area also has a high resolution. 

{\em Test 3.} The initial date is given as
\begin{align*}
	(\rho_0, u_0 ,v_0)=
	\left\{
	\begin{aligned}
		&(\rho_1, u_1, v_1) = (0.5287, 0.5, 0.5), \quad (x,y)\in [1,2]\times [1,2]\\
		&(\rho_2, u_2, v_2) = (0.7, 0.2, 0.5), \quad (x,y)\in [0,1]\times [1,2]\\
		&(\rho_3, u_3, v_3) = (0.4, 0.5, 0.5), \quad (x,y)\in [0,1]\times [0,1]\\
		&(\rho_4, u_4, v_4) = (0.7, 0.5, 0.2), \quad (x,y)\in [1,2]\times [0,1].
	\end{aligned}
	\right.
\end{align*}

	\begin{figure}[h!]
	\centering
	\vspace{-0.5cm}
		\includegraphics[width=\linewidth]{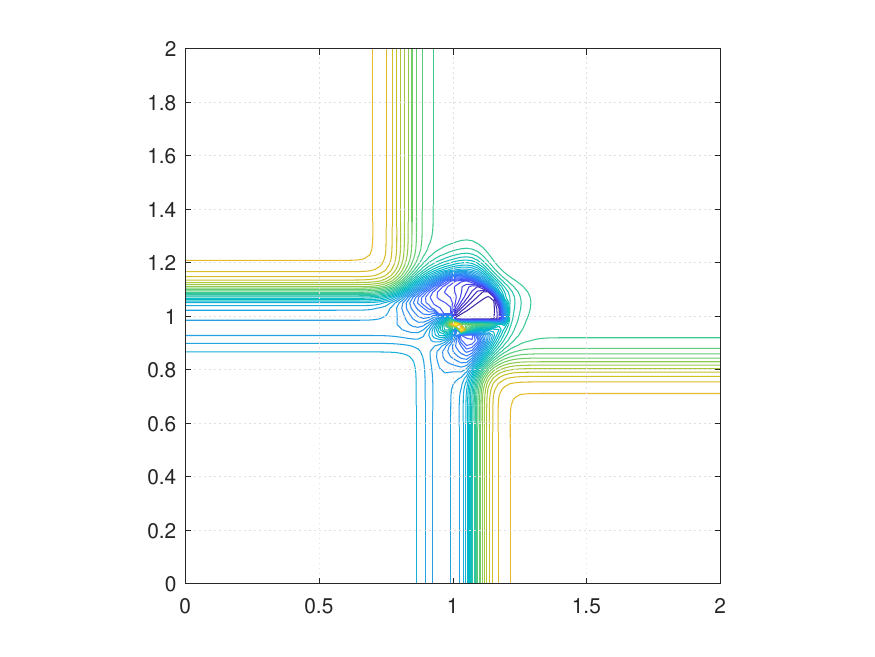}
		\caption{Density contours (HLL scheme)}
		\label{RJ_hll}
	\vspace{-0.1cm}
\end{figure}

\begin{figure}[h!]
	\centering
	\vspace{-0.5cm}
		\includegraphics[width=\linewidth]{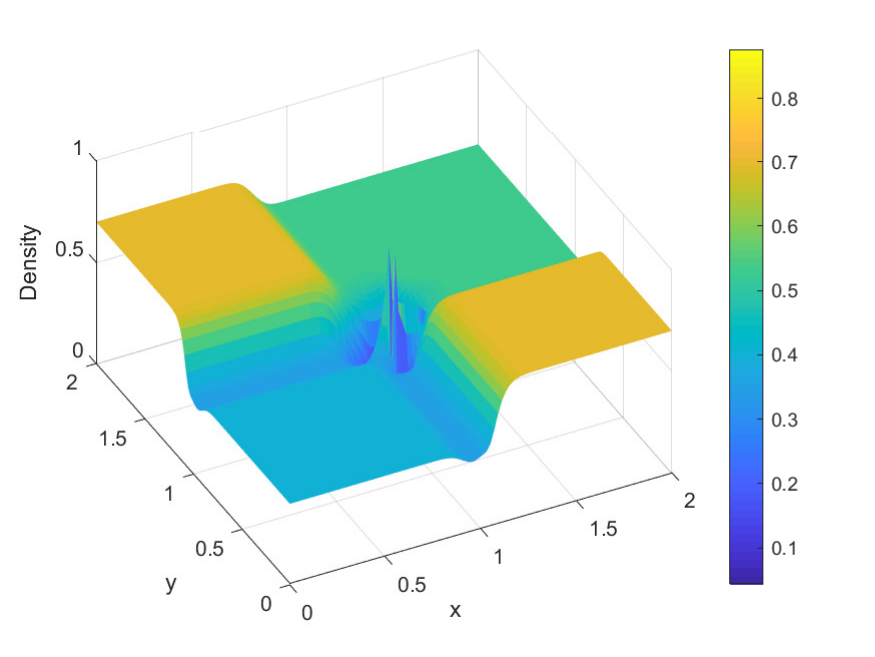}
		\caption{3D density field (HLL scheme)}
		\label{RJ_hll_3d}
	\vspace{-0.1cm}
\end{figure}

	\begin{figure}[h!]
	\centering
	\vspace{-0.5cm}
		\includegraphics[width=\linewidth]{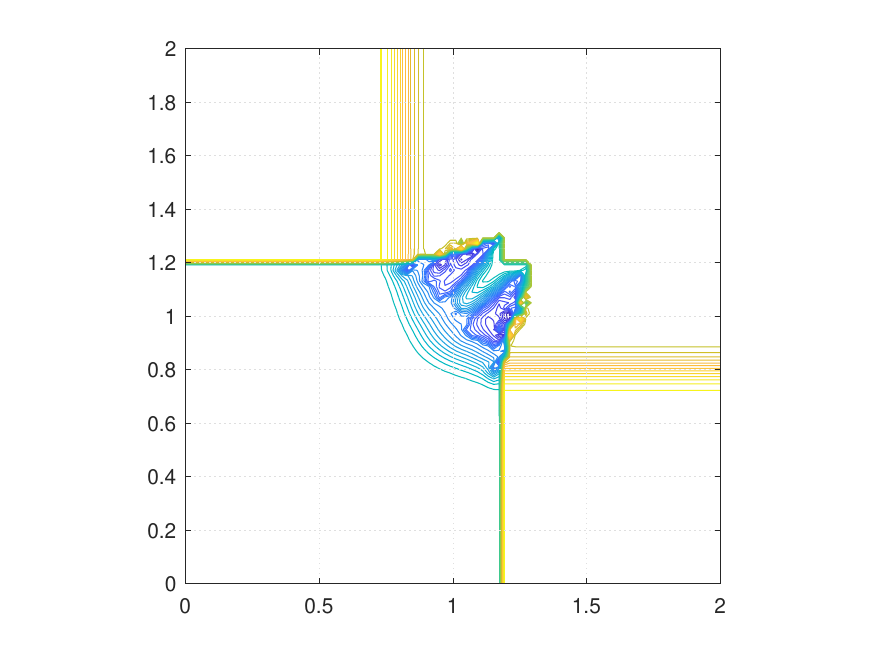}
		\caption{Density contours (Godunov-Glimm scheme)}
		\label{RJ_gg}
	\vspace{-0.1cm}
\end{figure}

\begin{figure}[h!]
	\centering
	\vspace{-0.5cm}
		\includegraphics[width=\linewidth]{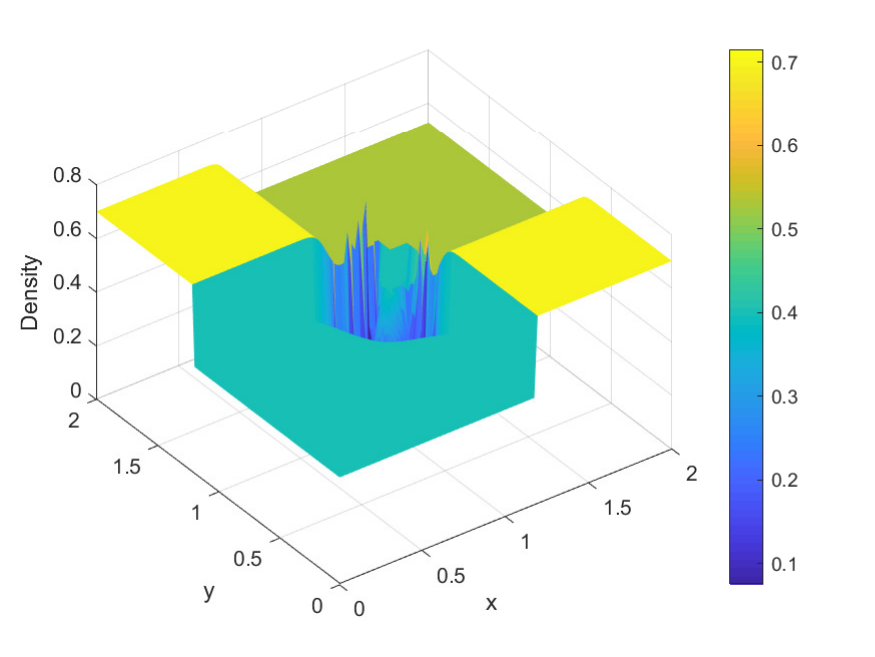}
		\caption{3D density field (Godunov-Glimm scheme)}
		\label{2R2J}
	\vspace{-0.1cm}
\end{figure}

The numerical solution is output at $t = 0.25\,\text{s}$, where the initial data consist of two rarefaction waves and two contact discontinuities, as shown in Fig.~\ref{2R2J}. Compared with the HLL scheme, the Godunov–Glimm method provides a better resolution, particularly in capturing the contact discontinuities and the interaction regions between different wave families. Again, it has been justified that the method can resolve the traffic dynamics more accurately,  and it offers  a faithful result of the underlying  traffic flow behavior.

\section{Conclusions}
In summary, this work proposes a refined traffic flow model by considering maximum density constraint and  maximum velocity constraint. 
From microscopic view, a quadratic term is introduced in the acceleration equation. The  formal macroscopic limit model indicates that the variant $\widetilde{u}p$ keeps constant along the traffic trajectories. The model’s fundamental diagram has been verified. The Riemann problem of the proposed model is analyzed and solved. Numerical simulations are carried on the hybrid Godunov-Glimm scheme, which confirm the accuracy and robustness of the model in reproducing key waves including shock waves, rarefaction waves, and stationary states. Beyond the one-dimensional case, we extend the model to two dimensions in order to account for more realistic traffic conditions, including multilane highways and urban road networks, demonstrating its adaptability to modern traffic systems.

%\bibliographystyle{siamplain}
%\bibliography{reference.bib}

\end{document}